\documentclass[a4paper,11pt]{article}
\usepackage{amssymb}
\usepackage{cite}

\newcommand{\s}{s}

\newcommand{\cL}{{\cal L}}

\newcommand{\cC}{\mathcal{C}}

\newcommand{\cK}{\mathcal{K}}
\newcommand{\cO}{\mathcal{O}}
\newcommand{\bX}{\mathbb{X}}

\newcommand{\R}{\mathbb{R}}

\newcommand{\N}{\mathbb{N}}

\newcommand{\cuad}{{\sqcap\kern-.68em\sqcup}}

\newcommand{\norm}[1]{\|#1\|}
\newcommand{\equ}[1]{(\ref{#1})}

\newtheorem{theorem}{Theorem}[section]
\newtheorem{proposition}{Proposition}[section]

\newtheorem{lemma}{Lemma}[section]
\newtheorem{corollary}{Corollary}[section]
\newtheorem{remark}{Remark}[section]

\begin{document}

\begin{center}{\bf   Isolated singularities for fractional Lane-Emden equations    \\[1.5mm]

  in  the Serrin's supercritical case

 }\bigskip

 {\small \begin{center}

 {\small
  Huyuan Chen\footnote{chenhuyuan@yeah.net}  \quad{and}\quad Feng Zhou\footnote{fzhou@math.ecnu.edu.cn}

\medskip
$^1$Department of Mathematics, Jiangxi Normal University,\\
 Nanchang, Jiangxi 330022, PR China \\[10pt]

 $^2$Center for PDEs, School of Mathematical Sciences, East China Normal University,\\
 Shanghai, 200241, PR China\\[19pt]
  }

\end{center}
}

\begin{abstract}
In this paper,  we  give a  classification of the isolated singularities of positive  solutions to the semilinear fractional elliptic  equations
$$(E) \qquad\qquad (-\Delta)^s u = |x|^{\theta} u^{p}\quad {\rm in}\ \ B_1\setminus\{0\},\qquad u= h\quad{\rm in}\ \ \R^N\setminus B_1,\qquad\qquad\quad   $$
where $s\in(0,1)$,  $\theta\in(-2s,0]$, $p>\frac{N+\theta}{N-2s}$, $B_1$ is the unit ball centered at the origin of $\R^N$ with $N>2s$.  $h$ is a nonnegative   function in $\R^N\setminus B_1$. Our   analysis of isolated singularities of $(E)$
is based on an integral upper bounds and the study of the Poisson problem with the fractional Hardy operators.  It is worth noting that our classification of isolated singularity holds in the
Sobolev super critical case $p>\frac{N+2s+2\theta}{N-2s}$ for $s\in(0,1]$ under suitable assumption of $h$.

\end{abstract}

  \end{center}
  \tableofcontents \vspace{1mm}
   \smallskip

  \noindent {\small {\bf AMS Subject Classifications}:  35R11; 35B40; 35J61.
  \smallskip

\vspace{2mm}

\setcounter{equation}{0}
\section{Introduction}
 Let  $s\in (0,1)$, $N>2s$ and  $B_r(y)$ be the ball of radius $r>0$ centered at $y$  in $\R^N$, $B_r:=B_r(0)$.
  Our concern of this paper is to classify  the isolated singular  positive solutions of the fractional Lane-Emden equation
\begin{equation}\label{eq 1.1}
\left\{
\begin{array}{lll}
(-\Delta)^s u=  |x|^{\theta} u^{p}  \quad &{\rm in}\ \, B_1 \setminus\{0\},
\\[2mm]
\qquad \ \ u=h&{\rm in}\ \,  \R^N\setminus B_1,
\end{array}\right.
 \end{equation}
where $\theta\in(-2s,0]$, $p>\frac{N+\theta}{N-2s}$,  $h$ is a nonnegative function in $\R^N\setminus B_1$ and
$(-\Delta)^s$ is the fractional Laplacian   defined by
$$
(-\Delta)^\s  u(x)= C_{N,\s}\lim_{\epsilon\to0^+} \int_{\R^N\setminus B_\epsilon }\frac{ u(x)-
u(x+z)}{|z|^{N+2\s}}  dz
$$
with
   $$C_{N,\s}=2^{2\s}\pi^{-\frac N2}\s\frac{\Gamma(\frac{N+2\s}2)}{\Gamma(1-\s)}$$
   being  the normalized constant, see \cite{NPV},    $\Gamma$ being the usual Gamma function.  It is known that  $(-\Delta)^\s  u(x)$ is well defined if $u$ is twice continuously differentiable in a neighborhood of $x$ and contained in the space   $L^1_{s}(\R^N):=L^1(\R^N,\frac{dx}{(1+|x|)^{N+2s}})$.
     \smallskip

When $s=1$ and $N\geq 3$, the isolated singularities of the Lane-Emden equation
\begin{equation}\label{eq 1.1-loc}
\left\{
\begin{array}{lll}
 -\Delta  u=  |x|^{\theta} u^{p}  \quad &{\rm in}\ \  B_1 \setminus\{0\},
\\[2mm]
\quad \ \, u=0&{\rm on}\ \,  \partial B_1
\end{array}\right.
 \end{equation}
 has been studied extensively in the last decades.  When  $\theta=0$ and $p\in(1,\, \frac{N}{N-2})$, Lions in \cite{L}
 classified the singular solution by building the connection with the weak solutions of
 $$-\Delta  u=   u^{p} +k\delta_0    $$
 and the related positive  solution $u_k$
 has  the asymptotic behavior
 $$\lim_{|x|\to0^+} u_k(x)|x|^{N-2}=c_N k $$
 for some $k\geq0$ and the normalized constant $c_N>0$ for $ -\Delta$.  In the particular case $k=0$, the solutions have removable singularity at the origin. In the
Serrin's critical and supercritical case $p\geq \frac{N}{N-2}$, however,  since it always has $k=0$, this method fails to classify the singularities.   When
  $$\theta\in(-2,2)\quad {\rm and}\quad \ p\in\Big(\frac{N+\theta}{N-2},\, \frac{N+2}{N-2}\Big)\setminus \Big\{\frac{N+2+2\theta}{N-2}\Big\},$$
   the singularity of solution for  (\ref{eq 1.1-loc}) was studied by Gidas-Spruck in \cite{GS} by using analytic techniques from \cite{A2}. In this case, the positive singular solutions have the singularity:
 $$u(x)= c_{p,\theta} |x|^{-\frac{2+\theta}{p-1}}(1+o(1))\quad{\rm as}\ \, |x|\to0^+$$
 with  the coefficient
\begin{equation}\label{eq 1.1-con}
c_{p,\theta}=\big(\frac{2+\theta}{p-1}(N-2-\frac{2+\theta}{p-1})\big)^{\frac{1}{p-1}}.\end{equation}
 When
 $$\theta\in(-2,2)\quad {\rm and}\quad \ p=\frac{N+\theta}{N-2},$$
  the author in \cite{A1,A2} gave a beautiful classification:  any positive solution $u$ of (\ref{eq 1.1-loc}) has removable singularity at origin or it has the asymptotic behavior
 $$u(x)=\Big(\frac{(N-2)^2}{2+\theta}\Big)^{\frac{N-2}{2+\theta}} |x|^{2-N}(-\ln|x|)^{\frac{2-N}{2+\theta}}(1+o(1))\quad{\rm as}\ \, |x|\to0^+.$$
 Moreover,   the existence of a singular solution is obtained by the phase plane analysis.  By the aid of  these solutions, Pacard in\cite{P} has constructed positive solutions with the prescribed singular set.
 Later on, Caffarelli-Gidas-Spruck in \cite{CGS} (also see \cite{KM}) classified the singular solutions of  (\ref{eq 1.1-loc}) when $\theta=0,\ p=\frac{N+2}{N-2}$ and its singular solutions satisfy
 $$u(x)=\varphi_{_D}(-\ln |x|) |x|^{-\frac{N-2}{2}}(1+o(|x|))\quad{\rm as}\ \, |x|\to0^+,$$
 where $\varphi_D$ may be the constant $c_{\frac{N+2}{N-2}}$ or be a periodic function.
In the  Sobolev super critical case $\theta=0$ and $p>\frac{N+2}{N-2}$, \cite{GNW,W} have studied the structure of positive radial solutions  of (\ref{eq 1.1-loc}), \cite{DGW} constructs
some non-radially symmetric solutions.  But  it is still open to give a full classification of isolated singularities for singular solutions.  We refer to \cite{BV,DDF,HLT,S0,S,V} for more  singularities for elliptic problem in various setting.

 Motivated by various applications and relationships to the theory of PDEs,  there has been an increasing interests in Dirichlet problems with nonlocal operators and the  prototype of nonlocal operator  is the fractional Laplacian.   The  problems with fully nonlinear nonlocal operators  are studied by \cite{CS0,CS1},  \cite{CS2} connects the fractional problem with the second order  degenerated problem in a half space with
  one dimensional higher,  and basic results about regularities for the fractional problem could see \cite{CS2,RS,Ls07} and  blowing up analysis for nonlocal problems could  see  \cite{CJ,CV1,CFQ}.
Motivated by \cite{L},  the isolated singularity of  (\ref{eq 1.1}) with $\theta=0$
was studied in \cite{CQ} for $p\in(1, \frac{N}{N-2s})$ via the connection with the distributional solutions  of
$$
(-\Delta)^s u=   u^{p} +k\delta_0 \quad  {\rm in}\ \, B_1 \setminus\{0\}, \qquad
  u=0 \ \ {\rm in}\ \,  \R^N\setminus B_1,
 $$
for which the positive singularity could be described by fundamental solution of fractional Laplacian.
As well as for the Laplacian, this method fails for the classification of isolated singularities for the Serrin's critical and supercritical case, i.e. $p\geq \frac{N}{N-2s}$.
 When $\theta=0$, $p\in(\frac{N}{N-2s},\frac{N+2s}{N-2s}]$,  \cite{JLX,CJ} build the platform of the isolated singularity for positive solution to
  \begin{equation}\label{eq 1.3}
\left\{
\begin{array}{lll}
 \,\quad\ \  {\rm div}(t^{1-2s}) \nabla U)=0  \quad &{\rm in}\ \, B_1\times (0,1)
\\[2mm]
\displaystyle  \lim_{t\to0^+}t^{1-2s}\partial_t U(x,t) =U^{p}(x,0)&  {\rm on}\ \,   B_1\setminus\{0\}
\end{array}\right.
 \end{equation}
 and  in the Sobolev's critical case $p=\frac{N+2s}{N-2s}$,  the behavior of non-removable singular solution of (\ref{eq 1.1})  can be stated as following
 $$c_2 |x|^{-\frac{N-2s}{2}}\leq u(x) \leq c_1|x|^{-\frac{N-2s}{2}}$$
  for some $c_1>c_2>0$,  and then \cite{YZ} gives a   description of singular solution near the singularity  for $p\in(\frac{N}{N-2s},\frac{N+2s}{N-2s})$ as
 $$c_4|x|^{-\frac{2s}{p-1}}\leq u(x)\leq c_3|x|^{-\frac{2s}{p-1}}$$
 for some $c_3>c_4>0$.  Also  its estimates singularity   for  $p=\frac{N}{N-2s}$ could see  \cite{WW}.      For the existence of isolated singular solutions,  \cite{AC,AD1} constructed a sequence of isolated singular solution of
 $$(-\Delta)^s u=u^p\quad{\rm in}\ \, \R^N\setminus\{0\}$$
 for   $p\in(\frac{N}{N-2s},\frac{N+2s}{N-2s})$, with fast decaying at infinity. We can also see \cite{AD,DG} for  the fractional Yamabe problem with isolated singularities in the Sobolev's critical case. Very recently, we provided in \cite{C} the classification of isolated singularities of (\ref{eq 1.1}) by
 a direct analysis method and showed the existence of a sequence of singular solutions
 of (\ref{eq 1.1}) with $\theta=0$ and $p=\frac{N}{N-2s}$.

Our aim in this paper is   to classify the positive singularity of fractional  problem (\ref{eq 1.1}) in the Serrin's supercritical case by a direct analysis method.  We first state our main results. To this end, we introduce some notations.
For $\tau\in(-N,2s)$, we denote
 \begin{equation}\label{cons 0}
  \cC_{s}(\tau) =2^{2\s}\frac{\Gamma(\frac{N+\tau}{2})\Gamma(\frac{2\s-\tau}{2})}{\Gamma(-\frac{\tau}{2})\Gamma(\frac{N-2\s+\tau}{2})}.
  \end{equation}
 More properties about $ \cC_{s}(\tau)$ can be found in the next section. In particular we notice that
 \begin{equation}\label{cons 1}
   \cC_{s}(\tau)>0\quad {\rm for}\ \ \tau\in (2s-N,0).
 \end{equation}
Let
   \begin{equation}\label{cons k}
  \cK_{p,\theta}:= \cC_{s}\Big(-\frac{2s+\theta}{p-1}  \Big)^{\frac{1}{p-1}},
 \end{equation}
then it's well defined since $-\frac{2s+\theta}{p-1} \in (2s-N, 0)$ for $p>\frac{N+\theta}{N-2s}$ and $\theta>-2s$. Here we use the notation $\cK_{p}=\cK_{p,0}$

\smallskip

 \vskip2mm
The classification of isolated singularities of (\ref{eq 1.1}) states as following.

\begin{theorem} \label{teo 0}
Assume that $s\in(0,1)$,  $N > 2s$,  $\theta=0$
    \begin{equation}\label{con s}
    p\in \Big(\frac{N }{N-2s}, \frac{N+2s }{N-2s}\Big)
 \end{equation}
 and $h\in   L^1_{s}(\R^N\setminus B_1)$ is  nonnegative.

  Let  $u_0$ be a  positive solution of (\ref{eq 1.1}), then $u_0$ has either removable singularity at the origin or
there holds
  \begin{eqnarray*}
\frac{\cK_{p}}{C_0} \leq  \liminf_{|x|\to0^+}u_0(x)  |x|  ^{\frac{2s }{p-1}}
  \leq  \cK_{p }   \leq \limsup_{|x|\to0^+}u_0(x) |x|  ^{\frac{2s }{p-1}}  \leq  C_0\cK_{p},
\end{eqnarray*}
where $C_0\geq 1$ is a constant from the Harnack inequality.
 \end{theorem}

 \begin{theorem} \label{teo 1}

 Assume that     $h$ is   radially symmetric and nonnegative  in $ \R^N\setminus B_1$,
 decreasing with respect to $|x|$  and $h(1)\geq0$,
\begin{equation}\label{con ss+0}
 \theta=0 \quad \ \&\quad \   p\geq  \frac{N+2s }{N-2s}
\end{equation}
 or
\begin{equation}\label{con ss-0}
\theta\in(-2s,0) \quad \ \&\quad \  p>\frac{N+\theta }{N-2s}.
\end{equation}

    Let  $u_0$ be a positive  solution of (\ref{eq 1.1}) such that
 \begin{equation}\label{con h}
  u_0(x)\geq  h(1)\quad{\rm for\ any}\ \,  x\in B_1,
  \end{equation}
 then   $u_0$ is  radially symmetric, decreasing with respect to $|x|$.
 Moreover,
 $u_0$ has either removable singularity at the origin or one has
  \begin{equation} \label{esta-1}
 \frac{\cK_{p,\theta}}{C_0} \leq  \liminf_{|x|\to0^+}u_0(x)  |x|  ^{\frac{2s+\theta}{p-1}}
  \leq  \cK_{p,\theta}   \leq \limsup_{|x|\to0^+}u_0(x) |x|  ^{\frac{2s+\theta}{p-1}}  \leq  C_0\cK_{p,\theta}.
 \end{equation}

\end{theorem}

 Our outline of the proofs is the following.   The first step is to get an bound
 \begin{equation}\label{esta-2-1}
 \limsup_{|x|\to0^+} u_0(x)|x|^{\frac{2s+\theta}{p-1}}<+\infty.
  \end{equation}
 To this end,  we build the following important upper estimate:
 \begin{proposition}\label{lm 3.1-1}
 Assume that $\theta\in(-2s,+\infty)$,  $p\geq \frac{N+\theta}{N-2s}$  and $u_0$ is a nonnegative classical solution of (\ref{eq 1.1}).
  Then we have
$
|x|^{\theta}  u_0^p\in L^1(B_1,(1-|x|)^s dx)
$
and  there exists a uniform constant $c_5>0$, independent of $r$ and $u_0$, such that
\begin{equation}\label{lp}
\int_{B_r} |x|^{\theta} u_0^p  dx\leq c_5 r^{N-\frac{\theta +2sp}{p-1}}, \quad \forall \ r\in(0,\frac{1}{2}).
\end{equation}
 \end{proposition}

With the help of Proposition \ref{lm 3.1-1},    we can  develop a direct blowing up technique for $\theta=0$, in which we don't have to transform
 the equation into an extension problem, see \cite{CJ,JX,YZ}, which doesn't require the boundedness in $L^1_s (\R^N, dx)$, thanks to the local property.  In our blowing up procedure,   if (\ref{esta-2-1}) fails,   the essential part is to keep the scaled singular solutions  uniformly bounded in $L^1_s (\R^N)$ for $p>\frac{N}{N-2s}$ by (\ref{lp}) with $\theta=0$ and then we pass to the limit to get a positive bounded solution for the limit equation and a contradiction arises from the  nonexistence of positive bounded solution for that limit equation.

  However, the blowing up technique fails for $\theta\not=0$, since the limit equation is
  $(-\Delta)^s u=0$ in $\R^N$ when $\theta<0$, which has positive solutions
  while for  $\theta>0$ there is no limit equation at all.

   Fortunately,  when $\theta\in(-2s,0)$, we can proceed a direct method of moving planes to obtain  the properties of symmetry and monotonicity in the $|x|$ direction under the assumption on $h$ of symmetry and monotonicity, motivated by \cite{FW}. As a consequence, these properties   improve (\ref{lp}) to our desired upper bound (\ref{esta-2}).

The second step is to get the estimate: for non-removable solution, there holds
\begin{equation}\label{mid in}
 \liminf_{|x|\to0^+}u_0(x)  |x|  ^{\frac{2s+\theta}{p-1}}
  \leq  \cK_{p,\theta}   \leq \limsup_{|x|\to0^+}u_0(x) |x|  ^{\frac{2s+\theta}{p-1}}\end{equation}
  for $p\in(\frac{N+\theta}{N-2s}, +\infty)$. 
These inequalities are motivated by the Poisson problem involving the fractional Hardy operators, where the coefficients of the Hardy potential is very sensitive to the blowing up rate of the related solutions.  Here we want to emphasize that the inequalities (\ref{mid in}) does not depends on the upper bounds and Harnack inequality.


Generally,   we can draw a conclusion for general $\theta>-2s$:
  \begin{theorem}\label{pr 1}
 Assume that $\theta\in(-2s,+\infty)$, $p > \frac{N+\theta}{N-2s}$ and $u_0$ is a nonnegative  classical solution of (\ref{eq 1.1}) verifying that
 there exist  $C_1\geq 1$ and $r_1\in(0,1)$ such that  for $x,y\in B_{r_1}\setminus\{0\}$
\begin{equation}\label{lp-2--}
u_0(x)\leq C_1 u_0(y)\quad {\rm for}\ \  \frac12\leq \frac{|x|}{|y|}\leq 2,
\end{equation}
then
  \begin{equation} \label{esta-2}
 \frac{\cK_{p,\theta}}{C_1} \leq  \liminf_{|x|\to0^+}u_0(x)  |x|  ^{\frac{2s+\theta}{p-1}}
  \leq  \cK_{p,\theta}   \leq \limsup_{|x|\to0^+}u_0(x) |x|  ^{\frac{2s+\theta}{p-1}}  \leq  C_1\cK_{p,\theta}.
 \end{equation}
 \end{theorem}

Note that (\ref{lp-2--}) is one type of  Harnack inequality, which, together with the pointwise upper bound, to obtain the lower bound  $\displaystyle \liminf_{|x|\to0^+}u_0(x) |x|  ^{\frac{2s+\theta}{p-1}} >0$ for the non-removable solution.  To make clear of our results,  we explain the relationship of integral upper bound (\ref{lp}), pointwise upper bound $\displaystyle \limsup_{|x|\to0^+}u_0(x) |x|  ^{\frac{2s+\theta}{p-1}} <+\infty$ and Harnack inequality  (\ref{lp-2--}): 
 \begin{enumerate}
\item[ $(i)$ ] {\it for $\theta=0$: integral upper bound \ + \ blow up analysis\quad  $\Rightarrow$ \quad pointwise upper bound\quad  $\Rightarrow$ \quad Harnack inequality;}

\item[$(ii)$]  {\it  for $\theta\in(-2s,0)$:   integral upper bound  \ + \ radial symmetry, increasing monotonicity \quad $\Rightarrow$\quad  pointwise upper bound\quad  $\Rightarrow$ \quad Harnack inequality; }

\item[ $(iii)$] {\it for $\theta\in(-2s,+\infty)$:  integral upper bound \ + \ Harnack inequality \quad $\Rightarrow$\quad  pointwise upper bound.}
    \end{enumerate}

\begin{remark}
 Under the assumption of $h$ in Theorem \ref{teo 1},    the method the moving plane works,  all such nonnegative solutions of (\ref{eq 1.1}) are radially symmetric and  decreasing with respect to $|x|$.

If     (\ref{eq 1.1}) has a  positive solution,  then
  \begin{equation}\label{eq 1.1-dual}
\left\{
\begin{array}{lll}
(-\Delta)^s u=  |x|^\theta u^p  \quad &{\rm in}\ \, B_1\setminus\{0\},
\\[2mm]
\qquad \ \ u=h&{\rm in}\ \,  \R^N\setminus B_1
\end{array}\right.
 \end{equation}
admits  the minimal nonnegative solution $u_{\rm min}$.

 Any solution $u\in C^s(\R^N)$ of (\ref{eq 1.1-dual})  verifies  that $\displaystyle \inf_{x\in B_1} u(x)\leq h(1)$, thans to $u(1)=h(1)$. A stronger version of (\ref{con h}) is
$$\inf_{x\in B_1} u_{\rm min}(x)= h(1).$$
Particularly,
$$u_{\rm min}\equiv0\quad {\rm if}\;\;  h\equiv0$$
and
$$u_{\rm min}>b\quad {\rm if}\quad h\equiv b\ \, {\rm for\  some \ } b>0.$$
Under these assumptions,  the isolated singularity of all nonnegative solution of (\ref{eq 1.1})  could  be classified.

  For the existence of minimal  solution  $u_{\rm min}$, it could be approached by the sequence of functions $\{v_m\}_{m\in\N}$, which are the solutions of
$$
\left\{
\begin{array}{lll}
(-\Delta)^s v_m=  |x|^\theta v_{m-1}^p  \quad &{\rm in}\ \, B_1,
\\[2mm]
\qquad \ \ v_m =h&{\rm in}\ \,  \R^N\setminus B_1
\end{array}\right.
$$
 with $v_0$ being the $s$-harmonic extension of $h$ in $B_1$,   if (\ref{eq 1.1}) has a nonnegative solution, which is an upper bound for $\{v_m\}$.
 \end{remark}

 Our methods could be extended to classify the singularity of prototype Lane-Emden equation
\begin{equation}\label{eq 1.1-lap}
\left\{
\begin{array}{lll}
 -\Delta  u=  |x|^{\theta} u^{p}  \quad &{\rm in}\ \, B_1 \setminus\{0\},
\\[2mm]
\quad\ \  u=h&{\rm in}\ \,  \partial B_1,
\end{array}\right.
 \end{equation}
 in the super critical case $p>\frac{N+2+2\theta}{N-2}$ and $\theta>-2$, where $h\geq0$.
 We have the following results:

  \begin{theorem} \label{teo 2}

 Assume that   $N\geq 3$, $\theta>-2$, $ p>\frac{N+2+2\theta }{N-2}$,
 $v_0$ is a positive  solution of (\ref{eq 1.1-lap})
 and $c_{p,\theta}$ is given in (\ref{eq 1.1-con}).


$(i)$ If  (\ref{lp-2--}) holds for $v_0$,  then $v_0$ has either removable singularity at the origin or one has
\begin{equation}\label{ss-s}
 \frac{c_{p,\theta}}{C_2} \leq  \liminf_{|x|\to0^+}v_0(x)  |x|  ^{\frac{2+\theta}{p-1}}
  \leq  c_{p,\theta}   \leq \limsup_{|x|\to0^+}v_0(x) |x|  ^{\frac{2+\theta}{p-1}}  \leq  C_2 c_{p,\theta},
 \end{equation}
 where $C_2\geq 1$ is a constant from the Harnack inequality.

$(ii)$  If $\theta\in(-2,0]$,   $h$ is a nonnegative constant on $\partial B_1$
 and
 \begin{equation}\label{con h-lap}
  v_0(x)\geq  h, \quad \, \forall\, x\in B_1,
  \end{equation}
 then   $v_0$ is  radially symmetric, decreasing with respect to $|x|$ and it is either removable   at the origin or has the singularity (\ref{ss-s}).

\end{theorem}

It is worth noting that the limit
 $$\lim_{|x|\to0^+}v_0(x) |x|  ^{\frac{2+\theta}{p-1}}=c_{p,\theta} $$
 is not always true in the Sobolev supercritical case. In fact,    when $\frac{N+1}{N-3}<p<p_c(N-1)$,  $p_c(N-1)$ being the  Joseph-Lundgren exponent,
  \cite{DGW} shows that
   $$
 -\Delta  u=    u^{p}  \quad  {\rm in}\ \, \R^N \setminus\{0\}$$
 has infinitely many non-radially symmetric solutions with the form that
 $$u(r,\omega)=c_p r^{-\frac{2}{p-1}}w_0(\omega),\quad (r,\omega)\in \R\times \partial B_1,$$
 where $w_0$ is non-constant solution in $ \partial B_1$. Together with our theorem \ref{teo 2}, we get a corollary that $c_{p,0}\in w_0[\partial B_1]$.

 \smallskip

  The remainder of this paper is organized as follows. In Section 2, we provide  preliminary estimates  of the fractional Hardy Poisson problem,  an essential integral upper bound and Harnack inequality.     In Section 3, we obtain the   upper    bounds of positive solutions of (\ref{eq 1.1})
  Section 4 is devoted to  isolated singularities of (\ref{eq 1.1})  and proofs of Theorem \ref{teo 0} and Theorem \ref{teo 1}.  
\setcounter{equation}{0}
\section{ Preliminary  }

Our analysis of isolated singularities for solutions with fractional laplacian is based on the fractional Hardy problem. In the sequel,   we do the $C^2$ extension for $h$ in $B_1$, till denoting it by $h$, such that
 $$h=0\quad {\rm in} \ B_{\frac34}.$$
 We denote
\begin{equation} \label{eqppoo}
  U_0(x)=-(-\Delta)^s  h(x)\quad {\rm for }\ x\in B_1,
  \end{equation}
  which is bounded locally in $B_1$.

\subsection{Poisson problem for fractional Hardy problem }

Recall that  for $\tau\in(-N,2s)$, we define
\begin{equation}
  \label{eq:fractional-power}
 \cC_\s(\tau) = 2^{2\s} \frac{\Gamma(\frac{N+\tau}{2})\Gamma(\frac{2\s-\tau}{2})}{\Gamma(-\frac{\tau}{2})\Gamma(\frac{N-2\s+\tau}{2})}= -\frac{c_{N,\s}}2 \int_{\R^N}\frac{|e_1+z|^{\tau}+|e_1-z|^\tau-2}{|z|^{N+2\s}}\,dz,
\end{equation}
where $e_1=(1,0,\cdots,0)\in\R^N$ and  there holds
\begin{equation}
  \label{eq:fractional-power-1}
 (-\Delta)^\s |\cdot|^\tau = \cC_\s(\tau)|\cdot|^{\tau-2\s} \quad  {\rm in}\ \   \R^N\setminus\{0\}.
\end{equation}
Note that $\cC_\s(\cdot)$
has two zeros point $\{0, 2s-N\}$, i.e.
$$\cC_\s(0)=\cC_\s(2s-N)=0.$$
 By \cite[Lemma 2.3]{CW}, we know that  the function $\cC_s$
is strictly concave and has a unique maximum at $\big\{\frac{2\s-N}{2}\big\}$ with the maximal value $2^{2\s}  \frac{\Gamma^2(\frac{N+2\s}4)}{\Gamma^2(\frac{N-2\s}{4})}.$
Moreover,
\begin{equation}
  \label{eq:c-symmetry}
  \cC_\s(\tau)= \cC_\s(2\s-N-\tau) \quad  {\rm for} \ \, \tau \in (-N,2\s)
\end{equation}
and
\begin{equation}
  \label{eq:c-asymptotics}
\lim_{\tau \to -N}\cC_\s(\tau) = \lim_{\tau \to 2 \s}\cC_\s(\tau)=-\infty.
\end{equation}
 Particularly, we have that
\begin{equation} \label{po-ne}
\cC_\s(\tau)>0\quad {\rm for}\ \tau\in (2s-N,0),\qquad \cC_\s(\tau)<0\quad {\rm for}\ \tau\in (-N,2s-N)\cup (0,2s).
\end{equation}

Now let
$$\cL^s_\mu u(x)=(-\Delta)^s u(x)+\frac{\mu}{|x|^{2s}}u(x), $$
where
$$\mu\geq \mu_0:=-2^{2\s}\frac{\Gamma^2(\frac{N+2\s}4)}{\Gamma^2(\frac{N-2\s}{4})}.$$
It is shown in \cite{CW} that   the linear equation
$$\mathcal{ L}_\mu^s u=0\quad{\rm in}\ \ \R^N\setminus \{0\}$$
 has two distinct radial solutions
 $$\Phi_{\mu}(x)=\left\{\arraycolsep=1pt
\begin{array}{lll}
 |x|^{\tau_-( \mu)}\quad
   &{\rm if}\;\; \mu>\mu_0\\[1.5mm]
 \phantom{   }
|x|^{-\frac{N-2s}{2}}\ln\left(\frac{1}{|x|}\right) \quad  &{\rm   if}\;\; \mu=\mu_0
 \end{array}
 \right.\quad  {\rm and}\;\;\  \Gamma_{ \mu}(x)=|x|^{\tau_+(\mu)},$$
where
$\tau_-( \mu)  \leq  \tau_+( \mu)$ verify that
 \begin{eqnarray*}
 &\tau_-(\mu)+\tau_+(\mu) =2\s-N \quad\  {\rm for\ all}\ \   \mu \ge \mu_0,\\[1.5mm]
&\tau_-( \mu_0)=\tau_+( \mu_0)=\frac{2\s-N}2,\quad\
\tau_-( 0)=2\s-N, \;\;\ \tau_+(0)=0,  \\[1.5mm]
&\displaystyle\lim_{\mu\to+\infty} \tau_-( \mu)=-N\quad {\rm and}\quad \lim_{\mu\to+\infty} \tau_+(\mu)=2\s.
 \end{eqnarray*}
In the sequel  of the paper and when there is no ambiguity,  we  use the notations $\tau_+=\tau_+(\mu)$, $\tau_-=\tau_-( \mu)$.
Moreover, we have that for $\mu\geq\mu_0$
$$\cC_s(\tau)+\mu>0\quad {\rm for}\ \ \tau\in(\tau_-,\tau_+)$$
and
$$\cC_s(\tau)+\mu<0\quad {\rm for}\ \ \tau\in(-N,\tau_-) \cup (\tau_+,2s).$$

Our analysis of singularities of positive solutions to (\ref{eq 1.1}) is based on the study of the Poisson problem
involving the fractional Hardy operator:
 \begin{equation}\label{eq 2.1fk}
\left\{
\begin{array}{lll}
\cL^s_\mu u=   f  \quad &{\rm in}\ \  B_1 \setminus\{0\},
\\[2mm]
\quad \   u=h&{\rm in}\ \   \R^N\setminus B_1.
\end{array}\right.
\end{equation}

The dual of the operator $\cL^\s_\mu$ is a weighted fractional Laplacian $(-\Delta)^\s_{\Gamma_\mu}$ given by
\begin{equation}\label{L}
(-\Delta)^\s_{\Gamma_\mu} v(x):=
C_{N,\s}\lim_{\epsilon\to0^+} \int_{\R^N\setminus B_\epsilon }\frac{v(x)-
v(z)}{|x-z|^{N+2\s}} \, \Gamma_\mu(z) dz.
\end{equation}
This expression is well defined for $x \in \R^N \setminus \{0\}$ if $v \in L^1(\R^N, \frac{\Gamma_\mu(x)}{1+|x|^{N+2\s}}dx)$ and if $v$ is twice continuously differentiable in a neighborhood of $x$.  From \cite[Proposition 3.1]{CW} there holds
\begin{equation}\label{1.3}
   | (-\Delta)^\s_{\Gamma_\mu}\xi(x)| \le  c_0 \min \{\Lambda _\mu(x),|x|^{-N-2\s}\}\quad  {\rm for} \ \, \xi\in C^2_c(\R^N),\ \,  \, x \in \R^N \setminus \{0\},
\end{equation}
where $c_0=c_0(\s,\mu,\xi)>0$ is a constant and
\begin{equation}\label{def-Lambda}
 \Lambda _\mu(x)= \left\{\arraycolsep=1pt
\begin{array}{lll}
  \displaystyle 1\quad
   &{\rm if}\ \,  \tau_+(\s,\mu)>2\s-1,\\[1mm]
   |x|^{1-2\s+\tau_+(\s,\mu)} \quad
   &{\rm if}\ \,  \tau_+(\s,\mu)<2\s-1,\\[1mm]
 \phantom{   }
 \displaystyle  1+(-\ln|x|)_+ \quad &{\rm   if } \ \,  \tau_+(\s,\mu) =2\s-1.
 \end{array}
 \right.
\end{equation}

 Then from \cite[Theorem 1.4]{CW} we have
\begin{theorem} \label{theorem-C}
 Let $\mu\ge\mu_0$ and $f\in C^\theta_{loc}(\bar B_1\setminus \{0\})$ for some $\theta\in(0,1)$.
\begin{enumerate}
\item[(i)] (Existence) If $f  \in L^1(B_1,\Gamma_\mu(x) dx)$, then for every $k\in\R$ there exists a solution $u_k \in L^1(B_1,\Lambda_\mu  dx)$ of   problem (\ref{eq 2.1fk})
satisfying the distributional identity
\begin{equation}
  \label{eq:distributional-k}
\int_{B_1}u_k   (-\Delta)^s_{\Gamma_\mu} \xi \,dx = \int_{B_1}f \xi\, \Gamma_\mu dx +c_{\s,\mu} k\xi(0) \quad  {\rm for\ all} \ \, \xi \in \cC^2_0(B_1).
\end{equation}
\item[(ii)]  (Existence and Uniqueness) If $f \in L^\infty(B_1, |x|^{\rho}dx)$ for some $\rho < 2\s - \tau_+(\s,\mu)$, then for every $k\in\R$ there exists a unique solution $u_k \in L^1(B_1,\Lambda_\mu  dx)$ of  problem (\ref{eq 2.1fk}) with the asymptotics
 \begin{equation}\label{beh 1}
 \lim_{|x| \to 0^+}\:\frac{u_k(x)}{\Phi_\mu(x)} = k.
 \end{equation}
Moreover, $u_k$ satisfies the distributional identity (\ref{eq:distributional-k}).
\item[(iii)] (Nonexistence) If $f$ is nonnegative and satisfies
\begin{equation}\label{f2}
 \int_{B_1} f\, \Gamma_\mu dx   =+\infty,
\end{equation}
then the problem
\begin{equation}\label{eq 1.1f}
 \arraycolsep=1pt\left\{
\begin{array}{lll}
 \displaystyle  \mathcal{L}_\mu^\s  u= f\quad
   {\rm in}\ \, B_1\setminus \{0\},\\[1.5mm]
 \phantom{  L_\mu \, }
 \displaystyle  u\ge 0\quad  {\rm   in}\ \,   \R^N\setminus B_1
 \end{array}\right.
\end{equation}
has no nonnegative distributional solution $u \in L^\infty_{loc}(\R^N \setminus \{0\}) \cap L^1_s (\R^N, dx)$.
\end{enumerate}
\end{theorem}

A direct consequence from Theorem \ref{theorem-C} with $\mu=0$ is the nonexistence of  solutions for
  \begin{equation}\label{eq 2.2-0}
  \left\{
\begin{array}{lll}
(-\Delta)^s u=   f  \quad &{\rm in}\ \, B_1 \setminus\{0\},
\\[2mm]
\qquad \  \   u=h&{\rm in}\ \,  \R^N\setminus B_1.
\end{array}\right.
\end{equation}

 \begin{corollary}\label{teo 2.1} Assume that the nonnegative function $f\in C^\beta_{loc}(\bar B_1\setminus\{0\})$ for some $\beta\in(0,1)$ and nonnegative function $h\in C^\alpha_{loc}(B_{2}\setminus B_1)\cap L^1(\R^N\setminus B_1, \frac{dx}{(1+|x|)^{N+2s}})$ with $\alpha\in(0,1)$.
 Then  problem (\ref{eq 2.2-0})   has no positive  solution,
  if
 $$\lim_{r\to0^+}\int_{B_1\setminus B_r}f(x)  dx=+\infty.$$
Particularly, the above assumption is filled if there holds
$$\liminf_{|x|\to0^+} f(x)|x|^{N} >0 .$$

\end{corollary}

Next we introduce  the  comparison principle of $\cL^s_\mu$.

\begin{lemma}\label{cr hp}
Assume that $\mu\geq\mu_0$, $O$ is a bounded $C^2$ domain containing the origin  and the functions $u_i\in C(\bar O\setminus\{0\})$ with $i=1,2$ verify in the classical sense
\begin{equation}\label{eq0 2.1}
 \arraycolsep=1pt\left\{
\begin{array}{lll}
 \displaystyle \cL^s_\mu u_1\leq  \cL^s_\mu u_2\quad
   &{\rm in}\quad  O\setminus \{0\},\\[2mm]
 \phantom{ \cL^s_\mu    }
 \displaystyle  u_1\leq u_2\quad  &{\rm   in}\ \ \R^N\setminus  O
 \end{array}\right.
\end{equation}
and satisfying
$$\limsup_{|x|\to0^+}u_1(x)|x|^{-\tau_-}\leq \liminf_{|x|\to0^+}u_2(x)|x|^{-\tau_-}.$$
Then
$$u_1\le u_2\quad{\rm in}\;\;  O\setminus \{0\}.$$

\end{lemma}
{\bf Proof.}    Let $u=u_1-u_2$ and then
$$\cL^s_\mu u \le 0\quad {\rm in}\ \, O\setminus\{0\}\quad {\rm and}\quad \limsup_{|x|\to0^+}u(x)\Phi_{s }^{-1}(x)\leq0,$$
 then for any $\epsilon>0$, there exists $r_\epsilon>0$ converging to $0$ as $\epsilon\to0$ such that
 $$u\le \epsilon \Phi_{\mu}\quad{\rm in}\;\; \overline{B_{r_\epsilon}(0)}\setminus\{0\},$$
 where $\Phi_{\mu}(x)=|x|^{\tau_-}$ is the fundamental solution of $\cL^s_\mu$.

We see that
$$u=0<\epsilon \Phi_{\mu} \quad{\rm in}\ \ \R^N\setminus O,$$
then we have
 $u\le \epsilon \Phi_{\mu} $ in $O\setminus\{0\} $ for any $\epsilon>0$, which implies that
$u\le 0$  in $O\setminus\{0\}$.
  \hfill$\Box$\medskip

   \begin{lemma}\label{teo 2.2} Assume that $\mu>\mu_0$,
     the nonnegative function $g\in C^\beta_{loc}(B_1\setminus\{0\})$ for some $\beta\in(0,1)$ and there exists $\tau\in\R $   such that
 $$ g(x)\geq |x|^{\tau-2s}\quad{\rm in}\ \,  B_{\frac12}\setminus\{0\}.$$
 Let  $u_g\in C(\bar B_1\setminus\{0\})$ verify
 \begin{equation}\label{eq 2.1-hom}
 \left\{\arraycolsep=1pt
\begin{array}{lll}
\mathcal{L}_\mu^\s  u_g\geq  g\quad\ \
   {\rm in}\ \,  B_1 \setminus \{0\},\\[2mm]
\quad \, u_g\geq  0\quad \ \  {\rm   on}\ \,  \R^N\setminus B_1.
  \end{array}
 \right.
 \end{equation}

We have the following
$(i)$ If $\tau\in(\tau_- ,\tau_+)$,
then there exists $c_6>0$  such that
$$u_g(x)\geq c_6|x|^\tau\quad{\rm in}\ \,  B_{\frac12}\setminus\{0\};$$

$(ii)$ If $\tau>\tau_+$,
then there exists $c_7>0$  such that
$$u_g(x)\geq c_7|x|^{\tau_+}\quad{\rm in}\ \,   B_{\frac12}\setminus\{0\}.$$

\end{lemma}
{\bf Proof. }   $(i)$ For $\tau\in(2s-N ,0)$, we have that
$$(-\Delta)^\s   |x|^{\tau}=\cC_s(\tau) |x|^{\tau-2s}\quad  {\rm in}\ \,  \R^N\setminus\{0\},$$
where $\cC_s(\tau)>0$ by (\ref{po-ne}).

We let
$$w_1(x)=|x|^{\tau}-1\quad  {\rm in}\ \  \R^N\setminus\{0\},$$
which, direct computation, verifies that
\begin{equation}\label{2.1-1}
\cL^s_\mu  w_1\leq (\cC_s(\tau)+\mu)|x|^{\tau-2}\quad{\rm in}\ B_1(0)\setminus\{0\},\qquad
w_1\leq 0\ \ {\rm in} \ \ \R^N\setminus B_1,
\end{equation}
where $\cC_s(\tau)+\mu>0$ for $\tau\in(\tau_- ,\tau_+)$ and $\mu>\mu_0$.

 Note that by the lower bound of $g$ there exists $t_0>0$ such that
 $$\cL^s_\mu  (t_0u_g)\geq t_0g(x)  \geq (\cC_s(\tau)+\mu)|x|^{\tau-2s}=\cL^s_\mu w_1\quad{\rm in}\ \, B_1\setminus\{0\}  $$
 and
 $$\liminf_{|x|\to0}u_g\Phi_\mu^{-1}(x)\geq0= \lim_{|x|\to0}w_1(x)\Phi_\mu^{-1}(x),\qquad t_0u_g\geq 0\geq w_1\quad{\rm in}\ \, \R^N\setminus B_1.$$
Then  by Lemma \ref{cr hp}, we have that
$$ t_0u_g\geq w_1\quad{\rm in}\ \, B_1\setminus\{0\}.$$

 $(ii)$ For $\tilde \tau\in\big( \tau_+,\min\{\tau, 2s\}\big)$ we set
  $$w_2(x)=|x|^{\tau_+}-|x|^{\tilde \tau}\quad  {\rm in}\ \  \R^N\setminus\{0\},$$
  then
$$\cL^s_\mu  w_2\leq -(\cC_s(\tilde \tau)+\mu)|x|^{\tilde \tau-2s}\leq -(\cC_s(\tilde \tau)+\mu)|x|^{ \tau-2s}\quad{\rm in}\ \, B_1\setminus\{0\}$$
and $$
w_2\leq 0\ \ {\rm in}\ \, \R^N\setminus B_1,
$$
where $ -(\cC_s(\tilde \tau)+\mu) >0$ for our choice of $\tilde \tau$ and $\tilde \tau-2s\geq   \tau-2s$.
Lemma \ref{cr hp} implies that
 $u_g\geq w_2$ in $B_1\setminus\{0\}$.
 \hfill$\Box$\medskip

\begin{lemma}\label{teo 2.3}  Let $\mu>\mu_0$, $g\in C^\beta_{loc}(B_1\setminus\{0\})$  with $\beta\in(0,1)$  be a   nonnegative function such that  there exists $\tau>\tau_- $  such that
 $$ g(x)\leq  |x|^{\tau-2s}\quad{\rm for}\ \, x\in B_1\setminus\{0\},$$
and $h$ be a nonnegative function in $C^\alpha (B_2\setminus B_1)\cap L^1_s(\R^N\setminus B_1)$ with $\alpha\in(0,1)$.

 Let  $u_g\in C(\bar B_1\setminus\{0\})$ verify
 \begin{equation}\label{eq 2.1-hom-}
 \left\{\arraycolsep=1pt
\begin{array}{lll}
\cL^s_\mu  u_g\leq g\quad\ \
   {\rm in}\ \   B_1 \setminus \{0\},\\[2mm]
\quad  \, u_g= h\quad \ \  {\rm   in}\ \  \R^N\setminus B_1,\\[2mm]
 \displaystyle\lim_{|x|\to0^+} u_g(x)|x|^{-\tau_-} =0.
  \end{array}
 \right.
 \end{equation}
 Then
$(i)$ if $\tau\in(\tau_- ,\tau_+)$
then there exists $c_8>0$  such that
$$u_g(x)\leq c_8|x|^\tau\quad{\rm for}\ \,x\in  B_1\setminus\{0\};$$
$(ii)$ if $\tau>\tau_+$
then there exists $c_{9}>0$  such that
$$u_g(x)\leq c_{9}|x|^{\tau_+} \quad{\rm for}\ \, x\in B_1\setminus\{0\}.$$

\end{lemma}
\noindent {\bf Proof. }  From (\ref{eqppoo}), a direct computation shows that for $x\in B_{\frac14}$
\begin{eqnarray}
 0<U_0(x) &=&  c_{N,s} \int_{\R^N\setminus B_{\frac12}}\frac{h(y)}{|x-y|^{N+2s}}dy\nonumber
 \\ &\leq &    c_{N,s}   \int_{\R^N\setminus B_{\frac12}}\frac{h(y)}{(|y|-\frac14)^{N+2s}}dy\leq c_{10}\|h\|_{L^1_s(\R^N)}
 \label{est 2.1}
 \end{eqnarray}
by using the fact that
$$|x-y|\geq |y|-\frac14\quad {\rm for}\ \, |x|\leq \frac14\ \, {\rm and}\ \, |y|\geq \frac12.$$

$(i)$ For $\tau\in(\tau_- ,\tau_+)$, we have that
$$\cL^s_\mu |x|^{\tau}=(\cC_s(\tau)+\mu) |x|^{\tau-2s}\quad  {\rm for}\ \, x\in \R^N\setminus\{0\},$$
where $\cC_s(\tau)+\mu>0$ and $\tau-2s<0$.

Note that  there exists $t_3>0$ such that
 $$\cL^s_\mu \big( t_3(u_g-h)\big)\leq t_3(g(x)+U_0) \leq  (\cC_s(\tau)+\mu) |x|^{\tau-2s} =\cL^s_\mu  |x|^{\tau} \quad  {\rm for}\ \, x\in B_\frac14\setminus\{0\},  $$
 $$\lim_{|x|\to0}\big(u_g(x)-h(x)\big)|x|^{-\tau_-}(x)=0 $$
  and
 $$t_3( u_g-h)_+\leq  |x|^{\tau}  \quad{\rm in}\ \ B_1\setminus B_\frac14\quad
{\rm and }\quad
 t_3( u_g-h) =0< |x|^{\tau}  \quad{\rm in}\ \ \R^N\setminus B_1.$$

Then  by Lemma \ref{cr hp}, we have that
$$ t_3\big(u_g(x)-h(x)\big)\leq  |x|^{\tau} \quad{\rm for}\ \, x\in B_1\setminus\{0\}.$$

$(ii)$    Take $\tilde\tau\in(\tau_+,2s)\cap(\tau_+,\tau]$ and
$$\tilde w_3(x)=  |x|^{\tau_+}-|x|^{\tilde\tau}+1\quad  {\rm in}\ \  B_1\setminus\{0\}\quad
{\rm and}\quad \tilde w_3(x)=1\quad{\rm in}\ \ \R^N \setminus B_1.$$
Direct computation shows that for $x\in B_{\frac{1}{2}}\setminus\{0\}$
\begin{eqnarray*}
\cL^s_\mu \tilde w_3(x)  &=&\cL^s_\mu \big(   |x|^{\tau_+}-|x|^{\tilde\tau}\big) +
c_{N,s}  \int_{\R^N\setminus B_1} \frac{ |y|^{\tau_+}-|y|^{\tilde\tau}}{|x-y|^{N+2s}}dy
\\[1.5mm]&\geq & - (\cC_{s}(\tilde \tau)+\mu) |x|^{\tilde \tau-2s}-c_{N,s}  \int_{\R^N\setminus B_1} \frac{|y|^{\tilde \tau} }{|x-y|^{N+2s}}dy
\\[1.5mm] &\geq &- (\cC_{s}(\tilde \tau)+\mu) |x|^{\tilde\tau-2s}-2^{N+2s}c_{N,s} \int_{\R^N\setminus B_1} \frac{|y|^{\tilde\tau} }{|y|^{N+2s}}dz,
\end{eqnarray*}
where $- (\cC_{s}(\tilde \tau)+\mu)>0$.

 Thus, there exists $r_0\in (0,\frac12]$ such that for $x\in B_{r_0}\setminus\{0\}$
 \begin{eqnarray*}
\cL^s_\mu  \tilde  w_3(x)  &\geq &-\frac12 (\cC_{s}(\tilde \tau)+\mu) |x|^{\tau-2s}.
\end{eqnarray*}

As a consequence,  there exists $t_4>0$ such that
 \begin{eqnarray*}\cL^s_\mu \big(t_4(u_g-h)\big)&\leq& t_4(g(x)+U_0)
 \\&\leq& -\frac12 (\cC_{s}(\tau)+\mu)  |x|^{\tilde\tau-2s}\leq \cL^s_\mu \tilde w_3(x) \quad{\rm in}\quad B_{r_0 }\setminus\{0\}
 \end{eqnarray*}
 and
 $$ t_4(u_g-h) \leq \tilde w_3\quad {\rm in}\ \, \R^N\setminus B_{r_0}.$$
Together with
 $$\lim_{|x|\to0^+}\big(t_4(u_g-h)(x)\big) |x|^{-\tau_-}=0,$$
 Lemma \ref{cr hp} is applied to obtain   that
$$ t_4(u_g-h) \leq  \tilde w_3 \quad{\rm in}\ \,    B_{r_0 }\setminus\{0\}.$$
Therefore,  we obtain that
$$u_g(x)\leq c_{9}|x|^{\tau_+} \quad{\rm for }\ \, x\in B_1\setminus\{0\}. $$
We complete the proof. \hfill$\Box$

 \subsection{Uniformly integral upper bounds}

 The integral upper bound (\ref{lp}) is  essential for our analysis of isolated singularities of  (\ref{eq 1.1})
 in the Serrin's supercritical case. \medskip

\noindent{\bf Proof of Proposition  \ref{lm 3.1-1}.  }   Recall $U_0$ in (\ref{eqppoo}),  which is uniformly  bounded in $B_{\frac12}$.
 Let  $w=u_0-h$ in $\R^N$, then we have that
 $$(-\Delta)^s w=|x|^{\theta} u_0^p+U_0\quad {\rm in}\ \, B_1,\qquad w=0\quad {\rm in}\ \R^N\setminus B_1.$$
 From Theorem \ref{theorem-C}, we have that $|x|^{\theta}  u_0^p+U_0\in L^1(B_1)$, so is $u_0$, thanks to the boundedness of $U_0$. Moreover, we have that
 $$\int_{B_1}w   (-\Delta)^\s \xi \,dx = \int_{B_1}(|x|^{\theta} u_0^p+U_0) \xi\,   dx \qquad  {\rm for\ all}\ \,\xi \in  C^2_0(B_1).   $$
Let $\xi_1$ be the positive Dirichlet eigenfunction of $(-\Delta)^s$ related to the first eigenvalue, subject to the zero Dirichlet boundary condition in $\R^N\setminus B_1$, i.e.
\begin{equation}\label{eq xi1}
 \left\{\arraycolsep=1pt
\begin{array}{lll}
(-\Delta)^\s  u=\lambda_{1}   u\quad\,
   &{\rm in}\ \,  B_1,\\[2mm]
\qquad \, u= 0\quad \ \  &{\rm   in}\ \,  \R^N\setminus B_1.
  \end{array}
 \right.
 \end{equation}
Then $\lambda_1>0$, $\xi_1$ is positive and  $\xi_1\in C^{2s}_{loc}(B_1)\cap C^s(\R^N)$. Moreover, by the standard argument we can obtain that $\xi_1\in C^2_0(B_1)$ and for some $c_{10}>1$ there holds
$$\frac1{c_{10}}\rho^s\leq  \xi_1\leq c_{10}\rho^s\quad {\rm in}\ \, B_1. $$
Using $\xi_1 $ as a test function, we have that
  \begin{eqnarray*}
  \int_{B_1} |x|^{\theta} u_0^p\xi_1\,  dx+\int_{B_1}U_0\xi_1\,   dx&=&\lambda_1\int_{B_1}w \xi_1\,dx
  \\&\leq &\Big(\int_{B_1}|x|^{\theta} w^p \xi_1\,dx\Big)^{\frac1p} \Big(\int_{B_1}  \xi_1|x|^{-\frac{\theta}{p-1}}\,dx\Big)^{1-\frac1p},
    \end{eqnarray*}
which implies
\begin{eqnarray*}
  \int_{B_1}|x|^{\theta}  w^p\rho^s\,dx \leq c_{11},
    \end{eqnarray*}
where $c_{11}>0$ depends on $h$ and  we used the fact that  $-\frac{\theta}{p-1}>-N,$
thanks to $p\geq \frac{N+\theta}{N-2s}$. Let
$$\xi_r(x)=  \xi_1(r^{-1}x)\quad {\rm for}\ \, x\in \R^N,$$
then $\xi_r$ is the solution of
\begin{equation}\label{eq xi}
 \left\{\arraycolsep=1pt
\begin{array}{lll}
(-\Delta)^\s  u=\lambda_{1} r^{-2s} u\quad\ \
   &{\rm in}\ \,  B_r,\\[2mm]
\qquad \ \ u= 0\quad \ \  &{\rm   in}\ \,  \R^N\setminus B_r.
  \end{array}
 \right.
 \end{equation}
  Let
 $$w=u_0-\tilde u_0,$$
 where
 $\tilde u_0$ is the $s$-harmonic extension of $u_0$ in $B_r$, i.e. the solution of
$$
 \left\{\arraycolsep=1pt
\begin{array}{lll}
(-\Delta)^\s  u= 0\quad\ \
   &{\rm in}\ \,  B_r,\\[2mm]
\qquad \ \ u= u_0\quad     &{\rm   in}\ \,  \R^N\setminus B_r.
  \end{array}
 \right.
$$
 Then we have that
$$
 \left\{\arraycolsep=1pt
\begin{array}{lll}
(-\Delta)^\s  w=|x|^{\theta} (w+\tilde u_0)^p\quad\ \
   &{\rm in}\ \,  B_r,\\[2mm]
\qquad \ \ w= 0\quad \ \  &{\rm   in}\ \,  \R^N\setminus B_r
  \end{array}
 \right.
$$
and
$$
  \int_{B_r}w(-\Delta)^s \zeta dx=\int_{\R^N}|x|^{\theta} (w+\tilde u_0)^p  \zeta dx,\quad \zeta\in \bX_r,
$$
  which implies that
$$
  \int_{B_r}u_0(-\Delta)^s \zeta dx=\int_{\R^N} |x|^{\theta} u_0^p\zeta dx+\int_{B_r} \tilde u_0 (-\Delta)^s \zeta dx,\quad \zeta\in \bX_r.
$$
Now we take $\zeta=\xi_r$, we have that
 $$
  \lambda_1 r^{-2s}\int_{B_r}u_0 \xi_r dx=\int_{\R^N} |x|^{\theta} u_0^p \xi_r dx+\int_{B_r} \tilde u_0 (-\Delta)^s \xi_r dx.
 $$

Again since $-\frac{\theta}{p-1}>-N$,  we have
  \begin{eqnarray*}
\int_{B_r} |x|^{\theta} u_0^p \xi_r dx &<&  \int_{B_r}  |x|^{\theta}  u_0^p\xi_r dx+\int_{B_r} \tilde u_0 (-\Delta)^s \xi_r dx
\\[1.5mm]&=&\lambda_1 r^{-2s}\int_{B_r}u_0\xi_r dx
\\[1.5mm]&\leq &\lambda_1 r^{-2s}\Big(\int_{B_r}  |x|^{\theta} u_0^p  \xi_r dx\Big)^{\frac1p}\Big(\int_{B_r}\xi_r|x|^{-\frac{\theta}{p-1}}  dx\Big)^{1-\frac1p}
\\[1.5mm]&\leq & \|\xi_1\|_{L^\infty(B_1)}^{1-\frac1p}\lambda_1 r^{(N-\frac{\theta}{p-1})(1-\frac1p)-2s}\Big(\int_{B_r}|x|^{\theta}u_0^p \xi_r dx\Big)^{\frac1p},
\end{eqnarray*}
which implies that
\begin{equation}\label{eq hext-3}
\int_{B_r} |x|^{\theta} u_0^p\xi_r dx < \|\xi_1\|_{L^\infty(B_1)} \lambda_1^{\frac{p}{p-1}} r^{N -\frac{\theta}{p-1}-\frac{ 2s p}{p-1}}.
  \end{equation}
Notice   that
$$\xi_r(x)=\xi_1(r^{-1}x)\geq \min_{z\in\bar B_{\frac12}}\xi_1(z)\quad {\rm in }\ B_{\frac r2},$$
then
  \begin{equation}\label{eq hext-3}
\int_{B_\frac r2}  |x|^{\theta} u_0^p  dx< c_{12}\int_{B_r}  |x|^{\theta} u_0^p\xi_r dx < c_{12}\|\xi_1\|_{L^\infty(B_1)} \lambda_1^{\frac{p}{p-1}} r^{N -\frac{\theta}{p-1}-\frac{ 2s p}{p-1}}.
  \end{equation}
 Replace $r$ by $\frac r2$, we obtain (\ref{lp}). \hfill$\Box$\medskip

With the help of Proposition  \ref{lm 3.1-1},  we have the following sharp upper bound
under some suitable restrictions of $u_0$.
\begin{corollary}\label{cr 3.1-1}
 Assume that $\theta>-2s$, $p\geq \frac{N+\theta}{N-2s}$ and $u_0$ is a nonnegative  classical solution of (\ref{eq 1.1}) verifying  the Harnack inequality:
 there exist  $C_0>1$ and $r_1\in(0,1)$ such that  for $x,y\in B_{r_1}\setminus\{0\}$
\begin{equation}\label{lp-2}
u_0(x)\leq C_0 u_0(y)\quad {\rm for}\ \ 1\leq \frac{|x|}{|y|}\leq 2,
\end{equation}
then  there exists a uniform $c_{13}>0$ independent of $r_1$ and $u$ such that
\begin{equation}\label{lp-1}
  u_0(x)\leq c_{13} |x|^{ -\frac{2s+\theta}{p-1}}\quad{\rm for}\ \, 0<|x|<1.
\end{equation}
 \end{corollary}
 \noindent{\bf Proof. } Take $|x|=r$, then (\ref{lp-2}) implies that  for any $y\in B_r\setminus B_{\frac{r}2}$, there holds that
 $$u_0(y)\geq \frac1{C_0}u_0(x). $$
From Proposition  \ref{lm 3.1-1}, we see that for any given $|x|=r$,
   \begin{eqnarray*}
C_0^{-p}\Big(1-(\frac12)^N\Big) |\partial B_1| \, u_0^p(x) r^{N-\theta}\leq  \int_{B_r\setminus B_\frac r2} |y|^\theta u_0^p(y)  dy<c_{14} r^{N-\frac{\theta}{p-1}- \frac{ 2s p}{p-1}},
 \end{eqnarray*}
 which implies that
 $$u_0 (x)\leq c_{13} |x|^{ -\frac{2s+\theta}{p-1}}. $$
We complete the proof. \hfill$\Box$

\subsection{Harnack inequality}
Our aim in this subsection is to obtain the Harnack inequality for singular solution
of (\ref{eq 1.1}).

\begin{proposition}\label{pr Harnack}
Assume that $\theta>-2s$, $p\geq \frac{N+\theta}{N-2s}$ and $u_0$ is a nonnegative solution of (\ref{eq 1.1}) such that
$$\limsup_{|x|\to0^+}u_0(x) |x|^{\frac{2s+\theta}{p-1}}<+\infty.$$

 Then there exists $C_0>0$ such that for all $r\in(0,\frac12)$
$$
 \sup_{x\in B_{2r}\setminus B_r} u_0(x)\leq C_0\Big(\inf_{x\in B_{2r}\setminus B_r} u_0(x)+\|u_0\|_{L^1_s(\R^N)}\Big).
$$
We assume more that $u_0$ is singular at the origin, then for all $r\in(0,\frac12)$
 \begin{equation}\label{harn in}
 \sup_{x\in B_{2r}\setminus B_r} u_0(x)\leq C_0 \inf_{x\in B_{2r}\setminus B_r} u_0(x).
 \end{equation}
 \end{proposition}
{\bf Proof. }  For fixed $r_0\in (\frac14,\frac12)$ and $x_0\in \R^N$ verifying $|x_0|=r_0$, from there exists $C>0$  independent of $u_0$  such that for any $t\in(0,\frac{r_0}4]$
 $$\sup_{x\in B_t(x_0)} u_0(x)<C.$$

 Let
 $$w(x)=u_0(x)(1-\eta_0(4x))-h(x) \ \ \ {\rm in} \ \ \R^N,$$
 where $\eta_0$ be a smooth function such that
  $\eta_0=1$ in $B_1(0)$ and $\eta_0=0$ in $\R^N\setminus B_2$,   we recall that $h=0$ in $B_{\frac34}$.  Then
 $$(-\Delta)^s w=|x|^\theta u_0^{p^*-1}w+(-\Delta)^s h+(-\Delta)^s (u_0\eta_0(4x))\quad {\rm for}\ \ x\in B_t(x_0),  $$
  where $(-\Delta)^s h>0$ in $B_1$.
 Note that
  $0<u_0^{p^*-1}\leq c_{15}$ for some $c_{15}>0$ independent of $u_0$ and
\begin{eqnarray*}
0&<&(-\Delta)^s h(x)+(-\Delta)^s (u_0(x)\eta_0(4x))
\\[1.5mm]&\leq& c_{15}\Big(\|h\|_{L^1_s(\R^N)}+\|u_0\|_{L^1_s(\R^N)}\Big)
\\[1.5mm]&\leq&  2c_{15}\|u_0\|_{L^1_s(\R^N)}.
\end{eqnarray*}
Then \cite[Theorem 1.1]{TX} (also see \cite[Theorem 11.1]{CS1}) implies that
$$\sup_{x\in B_t(x_0)} u_0(x)\leq C_3\Big(\inf_{x\in B_t(x_0)} u_0(x)+\|u_0\|_{L^1_s(\R^N )}\Big),$$
which infers
 \begin{equation}\label{harn in00}
 \sup_{x\in B_{2r_0}\setminus B_{r_0}} u_0(x)\leq C_3\Big(\inf_{x\in B_{2r_0}\setminus B_{r_0}} u_0(x)+\|u_0\|_{L^1_s(\R^N)}\Big)
\end{equation}
 by finite covering argument, the scaling property and the upper bound of $u_0$.

 Now we do the scaling:
 $$u_t(x)=t^{\frac{2s+\theta}{p-1}} u_0(tx) $$
 for $r\in(0,1]$.
 Then $u_t$ also verifies (\ref{eq 1.1}) and from the upper bound, we have
 $$r^{2s} |x|^{\theta} u_r(x)^{p-1}\leq C \quad {\rm for}\ \  r_0  t<|x|<2r_0t,$$
 where $C$ is dependent of $t$.

 It follows by (\ref{harn in00}) that
 $$
 \sup_{x\in B_{2r }\setminus B_{r }} u_0(x)\leq C_3\Big(\inf_{x\in B_{2r}\setminus B_{r }} u_0(x)+\|u_0\|_{L^1_s(\R^N)}\Big),
$$
where $r=r_0t$.
 Thanks to
$$\displaystyle \lim_{|x|\to0^+}  u_0(x)=+\infty,$$
we obtain   (\ref{harn in}).  We complete the proof.  \hfill$\Box$

\setcounter{equation}{0}
\section{Upper bounds }

\subsection{ The case: $\theta=0$}
 \begin{proposition}\label{pr 2.1}
 Assume that $\theta=0$ and
 $$  p\in\big(\frac{N }{N-2s},\frac{N+2s}{N-2s}\big). $$
 Let  $u_0$ be a nonnegative solution of (\ref{eq 1.1}), then
  there exists $c_{16}>0$ such that
 \begin{equation}\label{ll 3.0}
 u_0(x)\leq c_{16} |x|^{-\frac{2s }{p-1}} ,\quad\forall\, x\in B_1\setminus\{0\}.
 \end{equation}

 \end{proposition}

 \noindent{\bf Proof. } The similar upper bound could see \cite{YZ}, where the fractional laplacian is defined by the extension to a local setting.  Here we provide the direct blow-up analysis in the nonlocal setting.

 Suppose by contradiction that there exists a sequence of points $\{x_k\}\subset B_\frac12\setminus\{0\}$  such that $|x_k|\to0^+$ as $k\to+\infty$ and
$$
|x_k|^{\frac{2s }{p-1}}u_0(x_k) \to+\infty\quad{\rm as}\ \, k\to+\infty.
$$
We can choose  $x_k$ again having the property that
\begin{equation}\label{ll 3.1}
|x_k|^{\frac{2s }{p-1}}u_0(x_k)=\max_{x\in B_1\setminus B_{|x_k|}} |x|^{\frac{2s }{p-1}}u_0(x)\to+\infty\quad{\rm as}\ \, k\to+\infty
\end{equation}
by the fact that the mapping
$r\mapsto \max_{x\in B_1 \setminus B_r} |x|^{\frac{2s}{p-1}}u_0(x) $
is decreasing.

 We denote
$$\phi_k(x):= \big(\frac{|x_k|}{2}-|x-x_k| \big)  ^{\frac{2s }{p-1}}  u_0(x)\quad{\rm for}\ \, |x-x_k|\leq \frac{|x_k|}{2}.$$
Let $\bar x_k$ be the maximum point of $\phi_k$ in $B_{\frac{|x_k|}{2}}(x_k)$,
that is,
$$\phi_k(\bar x_k)=\max_{|x-x_k|\leq \frac{|x_k|}{2}} \phi_k(x).$$
Set
$$\nu_k=\frac12\big(\frac{|x_k|}{2}-|\bar x_k-x_k|\big),$$
then $0<2\nu_k<\frac{|x_k|}{2}$ and
$$  \frac{|x_k|}{2}-|x-x_k|\geq \nu_k\ \ {\rm for}\ \, |x-\bar x_k|\leq \nu_k.$$
By the definition of $\phi_k$,  for any $|x-\bar x_k|\leq \nu_k$,
$$(2\nu_k)^{\frac{2s }{p-1}}u_0(\bar x_k)=\phi_k(\bar x_k)\geq \phi_k(x_k)\geq \nu_k^{\frac{2s }{p-1}} u_0(x), $$
which implies that
$$2 ^{\frac{2s }{p-1}}u_0(\bar x_k)\geq u_0(x)\quad {\rm for\ any}\ |x-\bar x_k|\leq \nu_k.$$
Moreover, we see that
\begin{eqnarray*}
|\bar x_k|^{\frac{2s }{p-1}} u_0(\bar x_k)\geq
(2\nu_k)^{\frac{2s }{p-1}} u_0(\bar x_k)&=&\phi_k(\bar x_k)\geq \phi_k(x_k)\\[1.5mm]&\geq & (\frac{|x_k|}{2})^{\frac{2s }{p-1}}u_0(x_k)\to+\infty\quad {\rm as}\ \,  k\to+\infty
\end{eqnarray*}
by the fact that $|\bar x_k|\geq \frac{|x_k|}{2}\geq 2\nu_k$.

Denote
$$W_k(y)=\frac{1}{m_k}u_0\big( m_k^{-\frac{p-1}{2s }} y-\bar x_k \big),\quad\forall\, y\in \Omega_k\setminus \{X_k\},$$
where
$$m_k=u_0(\bar x_k),$$
 $$\Omega_k:=\Big\{y\in \R^N: \,  m_k^{-\frac{p-1}{2s }} y-\bar x_k\in B_1\Big\}$$ and
$$X_k= m_k^{\frac{p-1}{2s }}   \bar x_k.$$

Note that
$$|X_k|=\big(u_0(\bar x_k) |\bar x_k|^{\frac{2s }{p-1}} \big)^{\frac{p-1}{2s }}\to+\infty\quad {\rm as}\ \,  k\to+\infty.$$
Thus,   we have that for $y\in \Omega_k\setminus \{X_k\}$
\begin{eqnarray*}
(-\Delta)^s W_k(y)&=&\frac{1}{m_k^{p }}(-\Delta)^s u_0\big(  m_k ^{-\frac{p-1}{2s }} y-\bar x_k\big)
\\[1.5mm]&=&\frac{1}{m_k^{p } }  u_0^{p }\big(  m_k ^{-\frac{p-1}{2s }} y-\bar x_k\big)
\\[1.5mm]&=&   W_k^{p }(y),
\end{eqnarray*}
i.e.
\begin{equation}\label{eq lim}
(-\Delta)^s W_k(y)= W_k^{p}(y)\quad {\rm for}\ \ x\in \Omega_k\setminus \{X_k\}.
\end{equation}

We claim that there is $c_{17}>0$ independent of $k$ such that
$$\|W_k\|_{L^1_s(\R^N)} \leq c_{17}$$
and for any $\epsilon>0$, there exists $k_1>0$ and $R>0$ such that
\begin{equation}\label{l1 es}
\int_{\R^N\setminus B_R(0)} W_k(y)(1+|y|)^{-N-2s} dy \leq \epsilon.
\end{equation}

In fact,  since   $|\bar x_k|\to0$, we see that
\begin{eqnarray*}
0&\leq&  \int_{\R^N\setminus \Omega_k} W_k(y)(1+|y|)^{-N-2s} dy
\\[1.5mm]&=&\frac{1}{m_k}\int_{\R^N\setminus \Omega_k } u_0\big(  m_k ^{ -\frac{p-1}{2s } } y-\bar x_k\big) (1+|y|)^{-N-2s} dy
\\[1.5mm]&=&m_k^{-p} \int_{\R^N\setminus B_1 }h(z) \big( m_k^{-\frac{p-1}{2s }}+|z+\bar x_k|\big)^{-N-2s}dz
\\[1.5mm]&\leq & m_k^{-p}    \int_{\R^N  }h(z)  (1+|z+\bar x_k|)^{-N-2s}dz
\\[1.5mm]&\leq & u_0(\bar x_k)^{-p}    \int_{\R^N  }h(z)  (1+|z|)^{-N-2s}dz
\\[1.5mm]&\to & 0\quad{\rm as}\ \, k\to+\infty.
 \end{eqnarray*}
 Taking
 $$r_k=  |\bar x_k| \, u_0(\bar x_k)^{\frac{p-1}{ 2s }}\to+\infty\quad{\rm as}\ \, k\to+\infty,$$ we obtain that
 \begin{eqnarray*}
0&\leq&  \int_{B_{r_k}(X_k)} W_k(y)(1+|y|)^{-N-2s} dy\\[1.5mm]&\leq &\frac{r_k^{-N-2s}}{m_k }\int_{B_{r_k}(X_k)} u_0\big(m_k^{-\frac{p-1}{2s } } y -\bar x_k\big) dy
\\[1.5mm]&=& \frac{r_k^{-N-2s}m_k^{\frac{p-1}{2s }N} }{m_k } \int_{  B_{  |\bar x_k| }  } u_0 (z)   dz
\\[1.5mm]&\leq & r_k^{-N-2s}m_k ^{\frac{p-1}{2s }N-1}  \Big(\int_{  B_{  |\bar x_k| }  }  u_0^{p } (z)   dz \Big)^{\frac1{p }}  \Big(\int_{  B_{  |\bar x_k| }}  dz\Big)^{1-\frac1{p}}
\\[1.5mm]&=& c_{18} r_k^{-N-2s}m_k ^{\frac{p-1}{2s}N-1}  |\bar x_k|^{ N\frac{ p-1 }{p}} \Big(\int_{  B_{  |\bar x_k| }  }  u_0^{p } (z)   dz \Big)^{\frac1{p}}
\\[1.5mm]&\leq & c_{18}   r_k^{-N-2s}m_k^{\frac{p-1}{2s }N-1}  |\bar x_k|^{\frac{pN}{p-1}  + ( N -\frac{ 2s p}{p-1})\frac{1}{p}}
\\[1.5mm]&= &c_{18}   r_k^{-\frac{2s }{p-1}-2s}
  \\[1.5mm]& \to& 0\quad{\rm as}\ \, k\to+\infty,
 \end{eqnarray*}
where we used  (\ref{lp}) and the fact that  $\frac{\theta}{p-1}<N$.

 Moreover, $W_k(y)\leq 1$ in $\Omega_k\setminus B_{r_k}$, then
 \begin{eqnarray*}
   \int_{\Omega_k\setminus B_{r_k}(X_k)} W_k(y)(1+|y|)^{-N-2s} dy&\leq & \int_{\Omega_k\setminus B_{r_k}(X_k)}(1+|y|)^{-N-2s} dy
\\&\leq &  \int_{  \R^N} (1+|y|)^{-N-2s} dy
 \end{eqnarray*}
 and
 \begin{eqnarray*}
  && \int_{\big(\Omega_k\setminus B_{r_k}(X_k) \big)\setminus B_R(0)} W_k(y)(1+|y|)^{-N-2s} dy\\&\leq & \int_{\R^N\setminus B_{R}(0)}(1+|y|)^{-N-2s} dy
\\&\leq &  c_{19} R^{-2s},
 \end{eqnarray*}
 where $c_{19}>0$ is independent of $k$.
  Now we conclude (\ref{l1 es}) and the claim is proved. \medskip

Note that  $0<W_k\leq 2^{\frac{2s}{p-1}}$ in $B_{\tilde r_k}$,
where
$$\tilde r_k=\nu_km_k^{\frac{p-1}{2s}}=(\nu_k^{\frac{2+\theta}{p-1}}m_k)^{\frac{p-1}{2s}}\to+\infty\quad{\rm as}\ \, k\to+\infty,  $$
then for any $R>0$, there exists $k_R$ for any $k\geq k_R$
\begin{eqnarray*}
\|W_k\|_{C^{2s+\alpha}(B_{_R })} &\leq& c_{20}\Big( \|W_k\|_{L^1_s(\R^N)}+\|W_k\|_{L^\infty(B_{_{2R}})}+ \|W_k^p\|_{L^\infty(B_{_{2R}})}\Big)\\[2mm]&\leq& c_{20}\Big( \|W_k\|_{L^1_s(\R^N)}+2^{1+\frac{2}{p-1}p}\Big),
\end{eqnarray*}
where $\alpha\in(0,s)$ and $c_{20}>0$.

 By the arbitrary of $R$,  up to subsequence, there exists a nonnegative function $ W_\infty\in L^\infty(\R^N)$ such that as $k\to+\infty$
 \begin{equation}\label{cve1}
 W_k\to W_\infty\quad\ {\rm in}\ \ C^{2s+\alpha'}_{loc}(\R^N) \ \
 \ {\rm and\ \, in }\ \ L^1_s(\R^N)
 \end{equation}
 for some $\alpha'\in(0,\alpha)$.

For any $x\in \R^N$, take $R> 4|x|$ and $R\leq \min\big\{r_k,\tilde r_k\big\}$ for $k$ large enough,  note that
\begin{eqnarray*}
\frac1{C_{N,s}}(-\Delta)^s W_k(x)&=&{\rm p.v.}\int_{B_R} \frac{W_k(x)-W_k(y)}{|x-y|^{N+2s}} dy
+W_k(x) \int_{\R^N\setminus B_R} \frac{1}{|x-y|^{N+2s}} dy
\\&&\ \ -\int_{\R^N\setminus B_R} \frac{W_k(y)}{|x-y|^{N+2s}} dy
\\[1.5mm]&=:&E_{1,k}(x)+E_{2,k}(x)-E_{3,k}(x), \ \, {\rm respectively.}
\end{eqnarray*}
For any $\epsilon>0$,  by convergence (\ref{cve1}), we have that
$$\Big|E_{1,k}(x)-{\rm p.v.}\int_{B_R} \frac{W_\infty(x)-W_\infty(y)}{|x-y|^{N+2s}} dy\Big|<\frac\epsilon3$$
and there exists $R_1>0$ such that $R>R_1$
$$\Big| E_{2,k}(x)-W_\infty(x) \int_{\R^N\setminus B_R} \frac{1}{|x-y|^{N+2s}} dy\Big|<\frac\epsilon3.$$
Furthermore,  by (\ref{l1 es}) we have that
\begin{eqnarray*}
0<E_{3,k}(x)  &\leq &  \int_{\R^N\setminus B_R} W_k(y)(1+|y|)^{-N-2s} dy<\frac\epsilon3
\end{eqnarray*}
for $k>0$ and $R$ large enough.

Now we conclude that
$$\lim_{k\to+\infty}(-\Delta)^s W_k(x)=(-\Delta)^s W_\infty(x)$$
and $W_\infty$ is a bounded  classical solution of
\begin{equation}\label{eq lim1}
(-\Delta)^s W_\infty=W_\infty^{p}\quad {\rm in}\ \,\R^N.
\end{equation}

Since $W_\infty(0)=1$, then $W_\infty>0$ in $\R^N$,  thanks to the nonnegative property of $W_\infty$.
By   \cite[Theorem 3 ]{CLO} or  \cite[Theorem 4.5]{CLO1}, problem (\ref{eq lim1})
has no bounded positive solution for $p\in(1,\frac{N+2s}{N-2s})$. \smallskip
 We complete the proof. \hfill$\Box$\medskip


  \subsection{ The case: $\theta\in(-2s,0]$}
 When $\theta\in (-2s,0]$, we would derive the upper bound by
   the integral upper bound and increasing monotonicity. We recall (\ref{con ss+0}) and (\ref{con ss-0}) , i.e.
 $$\qquad\theta=0 \quad\&\quad  p\geq \frac{N+2s}{N-2s}$$
 and
 $$\theta\in(-2s,0) \quad\&\quad p>\frac{N+\theta}{N-2s}.$$
  The upper bound states as following.
   \begin{proposition}\label{pr 2.1-sup}
Assume that  (\ref{con ss+0}) or (\ref{con ss-0})  hold for $\theta$ and $p$,    function  $h$ is   radially symmetric  in $ \R^N\setminus B_1$ and
 decreasing with respect to $|x|$  and $h(1)\geq0$.

  Let  $u_0$ be a nonnegative solution of (\ref{eq 1.1}) satisfying (\ref{con h}),
    then
  there exists $c_{22}>0$ such that
 \begin{equation}\label{ll 3.0-sup}
0< u_0(x)\leq c_{22} |x|^{-\frac{2s+\theta}{p-1}} ,\quad\forall\, x\in B_1\setminus\{0\}.
 \end{equation}
 \end{proposition}

 In order to prove the upper bound, we first show the radial symmetry of solution of (\ref{eq 1.1}).
To this end, we use the moving plane method to show the radial symmetry and monotonicity of positive solutions to equation (\ref{eq 1.1}).

We  let
$$u=u_0-  h(1)\quad{\rm in}\ \ \R^N,$$
then $u\geq0$ in $B_1$ by (\ref{con h}) and it verifies that
\begin{equation}\label{eq 1.1-000}
\left\{
\begin{array}{lll}
(-\Delta)^s u=  |x|^{\theta} \big(u+ h(1)\big)^{p}  \quad &{\rm in}\ \, B_1 \setminus\{0\},
\\[2mm]
\qquad \ \ \, u=h-h(1)&{\rm in}\ \,  \R^N\setminus B_1.
\end{array}\right.
 \end{equation}

 \begin{proposition}\label{pr 2.1-sym}
Let  (\ref{con ss+0}), (\ref{con ss-0})  hold for $\theta$, $p$ respectively,    $h$ be   radially symmetric  in $ \R^N\setminus B_1$ and   decreasing with respect to $|x|$  and $h(1)\geq0$.

  Let  $u_0$ be a nonnegative solution of (\ref{eq 1.1}) satisfying (\ref{con h}),
then   $u=u_0-h(1)$ is a positive solution of   (\ref{eq 1.1-000}), $u$ is radially symmetric and decreasing with respect to
    $|x|$.
 \end{proposition}

 Our method of moving planes is motivated by \cite{FW} to deal with the solution $u$ has possibly singular at the origin.  For the methods of the moving planes for integral equations, we refer to \cite{CLO1,CLO}.
For singular solutions, we will use the direct moving plane method developing from  \cite{FW} and we need the following variant Maximum Principle for small domain.
\begin{lemma}\label{pr 2-a}\cite[Corollary 2.1]{FW} Let $O$  be an
open and bounded subset of $\R^N$. Suppose that
$\varphi\in L^\infty(O)$ and
$w\in L^\infty(\R^N)\cap C(\bar O)$ is a classical solution of
\begin{equation}\label{eq a1-a}
\left\{ \arraycolsep=1pt
\begin{array}{lll}
 (-\Delta)^s  w(x)\ge \varphi(x)w(x),\ \ \ \ &
x\in O,\\[2mm]
\qquad\ \  w(x)\ge 0,& x\in  O^c.
\end{array}
\right.
\end{equation}
Then there is $\delta>0$ such that whenever $|O^-|\le \delta$,
function $w$ has to be  non-negative in $O$, where $O^-=\{x\in O\ | \ w(x)<0\}$.
 \end{lemma}

  For simplicity, we denote
\begin{equation}\label{d1}\Sigma_\lambda=\{x=(x_1,x')\in \mathcal{O}_1\  |\  x_1>\lambda\},\end{equation}
\begin{equation}\label{d3}u_\lambda(x)=u(x_\lambda) \quad\mbox{and}\quad  w_\lambda(x)=u_\lambda(x)-u(x),\end{equation}
where $\lambda\in (0,1)$ and $x_\lambda=(2\lambda-x_1,x')$  for
$x=(x_1,x')\in\R^N$ and  $\mathcal{O}_1=B_1 \setminus\{0\}$.  For any subset $A$ of $\R^N$, we write
$A_\lambda=\{x_\lambda:\, x\in A\}$.

The essential estimate in the procedure is to show that for any $\lambda\in(0,1)$
$$w_\lambda\geq0 \quad{\rm in}\ \, \Sigma_\lambda.$$
On the contrary, we suppose that $\Sigma_\lambda^-=\{x\in \Sigma_\lambda\ |\ w_\lambda(x)<0\}\not=\emptyset$ for $\lambda\in(0,1)$.
Let us define
\begin{equation}\label{eq c}
w_\lambda^+(x)=\left\{ \arraycolsep=1pt
\begin{array}{lll}
w_\lambda(x),\ \ \ \ &
x\in \Sigma_\lambda^-,\\[1mm]
0,& x\in \R^N\setminus \Sigma_\lambda^-
\end{array}
\right.
\end{equation}
and
\begin{equation}\label{eq 4.01}
w_\lambda^-(x)=\left\{ \arraycolsep=1pt
\begin{array}{lll}
0,\ \ \ \ &
x\in \Sigma_\lambda^-,\\[1mm]
w_\lambda(x),\ \ \ \ & x\in \R^N\setminus \Sigma_\lambda^-.
\end{array}
\right.
\end{equation}
Hence, $w_\lambda^+(x)=w_\lambda(x)-w_\lambda^-(x)$ for all $x\in\R^N.$ It is obvious that
 $(2\lambda,0,\cdots,0)\not\in \Sigma_\lambda^-$, since $\displaystyle \lim_{|x|\to0^+}u(x)=+\infty$.

\begin{lemma}\label{lm 13-08-1}
Assume that   $\Sigma_\lambda^-\not=\emptyset$  for  $ \lambda\in(0,1)$, then
\begin{equation}\label{claim1}
 (-\Delta)^s w_\lambda^-(x)\le0,\ \ \ \ \forall\, x\in \Sigma_\lambda^-.
 \end{equation}

\end{lemma}
{\bf Proof.}  By direct computation, for
 $x\in \Sigma_\lambda^-$, we have
 \begin{eqnarray*}
 \frac1{C_{N,s}}(-\Delta)^s w_\lambda^-(x)
&=&\int_{\R^N}\frac{w_\lambda^-(x)-w_\lambda^-(z)}{|x-z|^{N+2s}}dz
=-\int_{\R^N\setminus\Sigma_\lambda^-}\frac{w_\lambda(z)}{|x-z|^{N+2s}}dz
\\&=&-\int_{(\mathcal{O}_1\setminus(\mathcal{O}_1)_\lambda) \cup ((\mathcal{O}_1)_\lambda\setminus \mathcal{O}_1)}\frac{w_\lambda(z)}{|x-z|^{N+2s}}dz
\\&&-\int_{(\Sigma_\lambda\setminus\Sigma_\lambda^-) \cup (\Sigma_\lambda\setminus\Sigma_\lambda^-)_\lambda}\frac{w_\lambda(z)}{|x-z|^{N+2s}}dz
-\int_{(\Sigma_\lambda^-)_\lambda}\frac{w_\lambda(z)}{|x-z|^{N+2s}}dz
\\&=&-I_1-I_2-I_3.
\end{eqnarray*}
We estimate these integrals separately. Since $u=0\ {\rm in} \ (\mathcal{O}_1)_\lambda\setminus \mathcal{O}_1$ and
$u_\lambda=0\ {\rm in} \ \mathcal{O}_1\setminus (\mathcal{O}_1)_\lambda$,
then
\begin{eqnarray*}
I_1
&=&\int_{(\mathcal{O}_1\setminus(\mathcal{O}_1)_\lambda) \cup ((\mathcal{O}_1)_\lambda\setminus \mathcal{O}_1)}\frac{w_\lambda(z)}{|x-z|^{N+2s}}dz
\\&=&\int_{(\mathcal{O}_1)_\lambda\setminus \mathcal{O}_1}\frac{u_\lambda(z)}{|x-z|^{N+2s}}dz-\int_{ \mathcal{O}_1\setminus(\mathcal{O}_1)_\lambda}\frac{u(z)}{|x-z|^{N+2s}}dz
\\&=&\int_{(\mathcal{O}_1)_\lambda\setminus \mathcal{O}_1}u_\lambda(z)\Big(\frac{1}{|x-z|^{N+2s}}-\frac{1}{|x-z_\lambda|^{N+2s}}\Big)dz\ge 0,
\end{eqnarray*}
since $u_\lambda\geq0$ and $|x-z_\lambda|>|x-z|$ for all  $x\in \Sigma_\lambda^-$ and $z\in (\mathcal{O}_1)_\lambda\setminus \mathcal{O}_1.$

In order to decide the sign of $I_2$ we observe that $w_\lambda(z_\lambda)=-w_\lambda(z)$ for any $z\in\R^N$. Then,
\begin{eqnarray*}
I_2
&=&\int_{(\Sigma_\lambda\setminus\Sigma_\lambda^-) \cup(\Sigma_\lambda\setminus\Sigma_\lambda^-)_\lambda}\frac{w_\lambda(z)}{|x-z|^{N+2s}}dz
\\&=&\int_{\Sigma_\lambda\setminus\Sigma_\lambda^-}\frac{w_\lambda(z)}{|x-z|^{N+2s}}dz
+\int_{\Sigma_\lambda\setminus\Sigma_\lambda^-}\frac{w_\lambda(z_\lambda)}{|x-z_\lambda|^{N+2s}}dz
\\&=&\int_{\Sigma_\lambda\setminus\Sigma_\lambda^-}w_\lambda(z)(\frac{1}{|x-z|^{N+2s}}-\frac{1}{|x-z_\lambda|^{N+2s}})dz
\\&\ge& 0,
\end{eqnarray*}
since  $w_\lambda\ge0$ in $\Sigma_\lambda\setminus\Sigma_\lambda^-$ and $|x-z_\lambda|>|x-z|$ for all $x\in \Sigma_\lambda^-$ and $z\in
\Sigma_\lambda\setminus\Sigma_\lambda^-.$

Finally, since
  $w_\lambda(z)<0$ for $z\in \Sigma_\lambda^-$, we deduce
\begin{eqnarray*}
I_3
&=&\int_{(\Sigma_\lambda^-)_\lambda}\frac{w_\lambda(z)}{|x-z|^{N+2s}}dz
=\int_{\Sigma_\lambda^-}\frac{w_\lambda(z_\lambda)}{|x-z_\lambda|^{N+2s}}dz
\\&=&-\int_{\Sigma_\lambda^-}\frac{w_\lambda(z)}{|x-z_\lambda|^{N+2s}}dz
\ge0.
\end{eqnarray*}
The proof is complete.\hfill$\Box$

\medskip

Now we are ready to prove Proposition \ref{pr 2.1-sym}.\medskip

\noindent{\bf Proof of  Proposition \ref{pr 2.1-sym}. }  For simplicity, let $h(1)=t\geq0$ and $u_0$ is nontrivial, i.e. $u_0\geq t$ $u_0\not\equiv t$ in $B_1$.

In order to show the radial symmetry and decreasing monotonicity in $|x|$,   we divide the proof into four steps.\smallskip

\noindent \emph{Step 1:} We prove that if $\lambda$ is close to $1$,  then  $w_\lambda>0$  in $\Sigma_\lambda$.  From the assumption  (\ref{con h}), we see that $w_\lambda>0$  in $\Sigma_\lambda$ for $\lambda=1$.

First we show that  $w_\lambda\geq0$  in $\Sigma_\lambda$, i.e. $\Sigma^-_\lambda$ is empty.
By contradiction, we assume that
$\Sigma^-_\lambda\not=\emptyset$. It must be in $B_1$
Now we apply  \equ{claim1} and linearity of the fractional Laplacian to   obtain that, for $ x\in\Sigma_\lambda^-,$
\begin{equation}\label{eq e}
(-\Delta)^s w_\lambda^+(x)\ge (-\Delta)^s w_\lambda(x)
=(-\Delta)^s u_\lambda(x)-(-\Delta)^s u(x).
\end{equation}
Combining    with (\ref{eq e}) and  (\ref{eq c}), for $x\in\Sigma_\lambda^-$, we have
\begin{eqnarray*}
(-\Delta)^s w_\lambda^+(x)  &\geq&  -|x_\lambda|^{\theta} (u_\lambda(x)+t)^{p}  +|x|^{\theta} (u (x) +t)^{p}
\\[2mm]&= &  -\varphi(x)w_\lambda^+(x),
\end{eqnarray*}
where
\begin{eqnarray*}
   \varphi(x)&=& \frac{  |x_\lambda|^{\theta} (u_\lambda(x)+t)^{p}-|x|^{\theta} (u(x)+t)^p }{u_\lambda(x)-u(x)}
   \\[2mm]&=&\frac{  |x_\lambda|^{\theta} -|x|^{\theta}   }{u_\lambda(x)-u(x)}(u_\lambda(x)+t)^{p} + |x|^{\theta} \frac{   (u_\lambda(x)+t)^{p}-  (u (x)+t)^{p} }{u_\lambda(x)-u(x)}
   \\[2mm]&< & |x|^{\theta}  \frac{   (u (x)+t)^{p}-  (u_\lambda(x)+t)^{p}}{(u(x)+t)-(u_\lambda(x)+t)}
   \\[2mm]&\leq &2^p  |x|^{\theta}  (u(x)+t)^{p-1},\quad\forall\, x\in\Sigma_\lambda^-
\end{eqnarray*}
thanks to $\theta\in(-2s,0)$.

For $x\in \Sigma_\lambda^-\subset \Sigma_\lambda\subset \R^N\setminus {B_\lambda }$, $u_\lambda(x)<u(x)$. Moreover,
there exists $M_\lambda>0$ such that
$$\norm{u}_{L^\infty(\R^N\setminus {B_\lambda})}\le M_\lambda,$$
   then
 there exists $c_{23}>0$ dependent of $\lambda$ such that
\begin{equation}\label{app 1}
\norm{\varphi}_{ L^\infty(\Sigma_\lambda^-)}\le c_{23}.
\end{equation}
Note that  $M_\lambda\to \infty$ as $\lambda\to 0$ if $\displaystyle \lim_{|x|\to0^+}u(x)=\infty$.\smallskip

Therefore, for  $x\in\Sigma_\lambda^-$ and then
$$
 (-\Delta)^s  w_\lambda^+(x)\ge -\varphi(x)w_\lambda^+(x), \quad
\forall\, x\in\Sigma_\lambda^-.
$$
Moreover, $w_\lambda^+=0$ in $(\Sigma_\lambda^-)^c$. Choosing $\lambda\in (0,1)$ close enough to $1$ we have $|\Sigma_\lambda^-|$ is small
and we apply Lemma \ref{pr 2-a} to obtain that
$$w_\lambda=w_\lambda^+\geq0\ \ \ \ \mbox{in} \ \ \Sigma_\lambda^-,$$
which is impossible. Thus,
$$w_\lambda\geq0\ \ \ \mbox{in}\ \ \Sigma_\lambda.$$

 If the function $\displaystyle\lim_{|x|\to0^+}u(x)=+\infty$, then   $w_\lambda$ is positive near the point $(2\lambda,0,\cdots,0)$ and then  $w_\lambda\not\equiv0$ in $\R^N$. If $u$ is bounded at the origin, then the solution $u$ has removable singularity at the origin.  In this case,
 $u=0$ on $\partial B_1$ then $w_\lambda\not\equiv 0$.

 Now we claim that for $0<\lambda<1$,    $w_\lambda>0$
 in $\Sigma_\lambda$.

Indeed, we assume on the contrary that there exists $x_0\in \Sigma_\lambda$ such that
$w_\lambda(x_0)=0,$ i.e.  $u_\lambda(x_0)=u(x_0)$. Then
\begin{eqnarray}\label{eq 7}
(-\Delta)^s w_\lambda(x_0)= (-\Delta)^s
u_\lambda(x_0)-(-\Delta)^s u(x_0)
=0.
\end{eqnarray}

On the other hand, let
$K_\lambda=\big\{(x_1,x')\in\R^N\ |\  x_1>\lambda\big\}$. Noting
$w_\lambda(z_\lambda)=-w_\lambda(z)$ for any $z\in\R^N$ and $w_\lambda(x_0)=0$, we deduce
\begin{eqnarray*}
(-\Delta)^s w_\lambda(x_0)
&=&-\int_{K_\lambda}\frac{w_\lambda(z)}{|x_0-z|^{N+2s}}dz-\int_{\R^N\setminus K_\lambda}\frac{w_\lambda(z)}{|x_0-z|^{N+2s}}dz
\\&=&-\int_{K_\lambda}\frac{w_\lambda(z)}{|x_0-z|^{N+2s}}dz-\int_{K_\lambda}\frac{w_\lambda(z_\lambda)}{|x_0-z_\lambda|^{N+2s}}dz
\\&=&-\int_{K_\lambda}w_\lambda(z)\Big(\frac{1}{|x_0-z|^{N+2s}}-\frac{1}{|x_0-z_\lambda|^{N+2s}}\Big)dz.
\end{eqnarray*}
The fact $|x_0-z_\lambda|>|x_0-z|$ for $z\in K_\lambda$ ,
$w_\lambda(z)\ge0$ and $w_\lambda(z)\not\equiv0$ in $K_\lambda$ yield
$$
(-\Delta)^s w_\lambda(x_0)<0,$$
which contradicts  (\ref{eq 7}), completing the proof of the claim.

\medskip

\noindent \emph{Step 2:} We prove $\lambda_0:=\inf\{\lambda\in(0,1)\ |\  w_\lambda>0\ \ \rm{in}\ \ \Sigma_\lambda\}=0$. If not, we set $\lambda_0>0$. Hence,
 $w_{\lambda_0}\geq0$ in $\Sigma_{\lambda_0}$ and
$w_{\lambda_0}\not\equiv0$ in $\Sigma_{\lambda_0}$. The claim in Step 1 implies
 $w_{\lambda_0}>0$ in $\Sigma_{\lambda_0}$.

Letting $\epsilon\in(0,\lambda_0/4)$ small enough, we now show
that $w_{\lambda_\epsilon }>0$ in $\Sigma_{\lambda_\epsilon}$ for
$\lambda_\epsilon=\lambda_0-\epsilon$.

To this end, we set $D_\mu=\{x\in\Sigma_\lambda\ | \ dist(x,\partial\Sigma_\lambda)\ge \mu\}$ for $\mu>0$ small. Since
$w_\lambda>0$ in $\Sigma_\lambda$ and $D_\mu$ is compact, there
exists $\mu_0>0$ such that $w_\lambda\ge \mu_0$ in $D_\mu$. By
the continuity of $w_\lambda(x)$, for $\epsilon>0$ small
enough and  $\lambda_\epsilon=\lambda-\epsilon,$ we have that
$w_{\lambda_\epsilon}(x)\ge0\ \ \rm{in}\ \ D_\mu$. Therefore,
$\Sigma_{\lambda_\epsilon}^-\subset
\Sigma_{\lambda_\epsilon}\setminus D_\mu$ and
$|\Sigma_{\lambda_\epsilon}^-|$ is small if $\epsilon$ and $\mu$ are small.
Using \equ{claim1} and proceeding as in {\it Step 1}, we have for all $x\in \Sigma_{\lambda_\epsilon}^-$
that \begin{eqnarray*}
(-\Delta)^s w_{\lambda_\epsilon}^+(x) &=&(-\Delta)^s
u_{\lambda_\epsilon}(x)-(-\Delta)^s u(x)-(-\Delta)^s
w_{\lambda_\epsilon}^-(x)
\\&\ge& (-\Delta)^s u_{\lambda_\epsilon}(x)-(-\Delta)^s u(x)
\\&\geq & -\varphi(x)w_{\lambda_\epsilon}^+(x).
\end{eqnarray*}
From (\ref{app 1})  $\varphi(x)$ is controlled by some constant dependent of $\lambda$.

Since $w_{\lambda_\epsilon}^+=0$ in
$(\Sigma_{\lambda_\epsilon}^-)^c$ and
$|\Sigma_{\lambda_\epsilon}^-|$ is small, for  $\epsilon$ and $\mu$
small,  Lemma \ref{pr 2-a}  implies that  $w_{\lambda_\epsilon}\ge0$ in
$\Sigma_{\lambda_\epsilon}$.  Combining with $\lambda_\epsilon>0$ and
$w_{\lambda_\epsilon}\not\equiv0$ in $\Sigma_{\lambda_\epsilon}$, we obtain
 $w_{\lambda_\epsilon}>0$
 in $\Sigma_{\lambda_\epsilon}$.   This contradiction arises from the definition of $\lambda_0$.

Therefore,  we have  that $\lambda_0=0$.

\medskip

\noindent \emph{Step 3:} By Step 2, we have $\lambda_0=0$, which
implies that $u(-x_1,x')\ge u(x_1,x')$ for $x_1\ge0.$
Using the same argument from the other side, we conclude that $u(-x_1,x')\le u(x_1,x')$ for $ x_1\ge0$ and then
$u(-x_1,x')= u(x_1,x')$ for $x_1\ge0.$ Repeating this procedure in all directions we see that $u$ is radially symmetric.

\medskip

Finally, we prove $u(r)$ is strictly decreasing in  $r\in (0,1)$. Let us consider  $0<x_1<\widetilde{x}_1<1$
and let  $\lambda=\frac{x_1+\widetilde{x}_1}{2}$. As proved above   we have
$$w_\lambda(x)>0\ \ \mbox{for}\ \ x\in\Sigma_{\lambda}.$$
Then
\begin{eqnarray*}
0<w_\lambda(\widetilde{x}_1,0,\cdots,0)
&=&u_\lambda(\widetilde{x}_1,0,\cdots,0)-u(\widetilde{x}_1,0,\cdots,0)
\\&=&u(x_1,0,\cdots,0)-u(\widetilde{x}_1,0,\cdots,0),
\end{eqnarray*}
i.e
$u(x_1,0,\cdots,0)>u(\widetilde{x}_1,0,\cdots,0).$
From the radial symmetry of $u$ and decreasing in the direction $\frac{x}{|x|}$,  we can conclude the monotonicity of $u$.
 \hfill$\Box$\medskip\medskip

\noindent {\bf Proof of Proposition \ref{pr 2.1-sup}. }  Note that
$$u=u_0- h(1)$$
and Proposition \ref{pr 2.1-sym} shows that   $u$ is radially symmetric and decreasing with respect to  $|x|$, so is $u_0$. Then (\ref{lp-2}) holds for $u_0$ and we apply  Corollary \ref{cr 3.1-1} to obtain that
$$u_0(x)\leq c_{23}|x|^{-\frac{2s+\theta}{p-1} }\quad{\rm for}\ \, x\in B_1\setminus\{0\}.$$
 We complete the proof. \hfill$\Box$

    \setcounter{equation}{0}
    \section{Classification of isolated singularity}
   \subsection{ Some important estimates }

   \begin{proposition} \label{teo 1-sub}
Assume that  $h\in   L^1_s(\R^N\setminus B_1)$,   $\theta\in(-2s,+\infty)$
and
  $$p\in \Big(\frac{N+\theta}{N-2s}, \frac{N+2s+2\theta}{N-2s}\Big).$$

Let $u_0$ be a  positive solution of (\ref{eq 1.1}) verifying
$$\limsup_{|x|\to0^+} u_0(x) =+\infty, $$
 then $u_0$ satisfies that
  \begin{eqnarray*}
   \liminf_{|x|\to0^+}u_0(x)  |x|  ^{\frac{2s+\theta}{p-1}}
  \leq  \cK_{p,\theta}   \leq \limsup_{|x|\to0^+}u_0(x) |x|  ^{\frac{2s+\theta}{p-1}},
\end{eqnarray*}
  where $\cK_{p,\theta}$ is given in (\ref{cons k}).
 \end{proposition}

\noindent{\bf Proof. }
We recall that
$$\cK_{p,\theta}  =  \cC_{s}\big(-\frac{2s+\theta}{p-1}\big)^{\frac{1}{p-1}}.$$
For $p\in\big(\frac{N+\theta}{N-2s}, \frac{N+2s+2\theta}{N-2s})$, there holds that
$$-\frac{2s+\theta}{p-1}\in(2s-N,\frac{2s-N}{2})$$
and
$$\tau_-(-\cK_{p,\theta}  ^{p-1})=-\frac{2s+\theta}{p-1}\quad {\rm and}\quad \tau_+(-\cK_{p,\theta}^{p-1})=2s-N+\frac{2s+\theta}{p-1}.$$

From Proposition \ref{pr 2.1}, we have that
$$\limsup_{|x|\to0^+} u_0(x)|x|^{\frac{2s+\theta}{p-1}}<+\infty.$$

  Set
\begin{equation}\label{1.2-400}
k= \liminf_{|x|\to0^+}u_0(x)  |x|^{\frac{2s+\theta}{p-1}}.
 \end{equation}

{\bf Part 1: } we claim
 $$  k\leq \cK_{p,\theta}. $$
In fact, if $$k\in(\cK_{p,\theta},+\infty],$$
 letting
$$\epsilon_0=\frac{k-\cK_{p,\theta} }{2\cK_{p,\theta}  },$$
then  there exists $r_1\in(0,1)$ such that
$$u_0(x)\geq \cK_{p,\theta}   \big(1+\epsilon_0\big)  |x|^{-\frac{2s+\theta}{p-1}} \quad {\rm for}\ \ 0<|x|<r_1, $$
then for $0<|x|<r_1$,
\begin{eqnarray*}
 u_0^{p-1}(x)&\geq& \cK_{p,\theta}^{p-1} \big(1+\epsilon_0\big)^{p-1} |x|^{-\theta-2s}>\cK_{p,\theta}  ^{p-1} \big(1+(p-1)\epsilon_0\big)  |x|^{-\theta-2s}.
 \end{eqnarray*}

 Therefore,     $u_0$ verifies that
\begin{equation}\label{eq 6.1-dd}
\left\{
\begin{array}{lll}
(-\Delta)^s u_0= \frac{\cK_{p,\theta}^{p-1}  +\frac{p-1}2\epsilon_0  }{|x|^{2s} } u_0+ f_1    \quad &{\rm in}\ \ B_{r_1} \setminus\{0\},
\\[2mm]
\qquad \ \   u_0\geq 0&{\rm in}\ \  \R^N\setminus B_{r_1},
\end{array}\right.
\end{equation}
where
$$\mu_p=-\Big(\cK_{p,\theta}^{p-1}  +\frac{p-1}2\epsilon_0\Big)<-\cK_p^{p-1}$$ and
\begin{eqnarray*}
f_1(x) &=&|x|^\theta \Big(u_0^{p-1}- \frac{\cK_{p,\theta}^{p-1}  +\frac{p-1}2\epsilon_0  }{|x|^{2s+\theta} } \Big) u_0(x)
\\& \geq& \frac{p-1}2 \epsilon_0\cK_{p,\theta}^p |x|^{-\frac{(2s+\theta)p}{p-1}}.
\end{eqnarray*}
Observe that
$$\tau_+(\mu_p)<\tau_+(-\cK_{p,\theta}^{p-1} )=2s-N+\frac{2s+\theta}{p-1}$$
that is
$$\tau_+(\mu_p)-\frac{(2s+\theta)p}{p-1}<-N, $$
which implies that
$$f_1\not\in L^1(B_{r_1},\Gamma_{\mu_p}dx)$$
and a contradiction arises by Theorem \ref{theorem-C}, from which problem (\ref{eq 6.1-dd}) has no such positive solution. Therefore we obtain  that
  $  k\leq \cK_p. $\smallskip

 {\bf Part 2:}   Set
\begin{equation}\label{1.2-40}
\kappa= \limsup_{|x|\to0^+}u_0(x)  |x|^{\frac{2s+\theta}{p-1}}
 \end{equation}
 and    we prove that
 $$  \kappa\geq \cK_{p,\theta}. $$
In fact, if
$$\kappa<\cK_{p,\theta}, $$ letting
$$\epsilon_1=\min\Big\{\frac{\cK_{p,\theta}-\kappa}{2\cK_{p,\theta}},\, \frac14\Big\},$$
then  there exists $r_2\in(0,1)$ such that
$$u_0(x)\leq \cK_{p,\theta} \big(1-\epsilon_1) |x| ^{-\frac{2s+\theta}{p-1}} \quad {\rm for}\ \ 0<|x|<r_2, $$
and for $0<|x|<r_2$
\begin{eqnarray*}
 u_0^{p-1}(x)&\leq& \cK_{p,\theta}^{p-1}\big(1-\epsilon_1)^{p-1}  |x| ^{-\theta-2s}.
 \end{eqnarray*}
 Let
 $$\tilde \mu_{p}=-\cK_{p,\theta}^{p-1}\big(1-\epsilon_1)^{p-1}>-\cK_{p,\theta}^{p-1},$$
 then
 $$\tau_-(\tilde \mu_{p})<-\frac{2s+\theta}{p-1}.$$

Then $u_0$ is a positive solution of
 \begin{equation}\label{eq 6.1-dd}
\left\{
\begin{array}{lll}
(-\Delta)^s u_0= \frac{\cK_{p,\theta}^{p-1}\big(1-\epsilon_1)^{p-1} }{|x|^{2s} } u_0+f_2    \quad &{\rm in}\ \, B_{r_2} \setminus\{0\},
\\[2mm]
\qquad \ \   u_0\geq 0&{\rm in}\ \,  \R^N\setminus B_{r_2},\\[2mm]
\displaystyle\lim_{|x|\to0^+ } u_0(x)|x|^{-\tau_-(\tilde \mu_{p})}=0,
\end{array}\right.
\end{equation}
where
\begin{eqnarray*}
f_2(x) =|x|^\theta \Big(u_0^{p-1}-  \frac{\cK_{p,\theta}^{p-1}\big(1-\epsilon_1)^{p-1} }{|x|^{2s+\theta} }\Big) u_0(x)
 \leq 0.
\end{eqnarray*}
Then Lemma \ref{teo 2.3}  with $\mu=\tilde \mu_p$  implies  that
$$\limsup_{|x|\to0^+}u_0(x)|x|^{-\tau_+(\tilde \mu_{p})}<+\infty, $$
where
$$\tau_+(\tilde \mu_{p})>\tau_+(-\cK_{p,\theta}^{p-1})=2s-N+\frac{2s+\theta}{p-1}>-\frac{2s+\theta}{p-1}.$$

Now  we take the value $\tau_0=\tau_+(\tilde \mu_{p})<0$,    then
 $$
(-\Delta)^s  u_0(x)  \leq c_{24}|x|^{\tau_0 p+\theta}    \quad {\rm in}\ \, B_{r_2} \setminus\{0\}
$$
and letting
$$\tau_1:=p\tau_0 +\theta+2s,$$
if $\tau_1> 0$,
 $$u_0(x)\leq d_1 \quad {\rm in}\ \, B_{r_2} \setminus\{0\},$$
 which contradicts  that $u_0$ is non-removable at the origin. Then we are done.

If $\tau_1=0$, applying Lemma \ref{teo 2.3} with $\mu=0$ to obtain that
$$u_0(x)\leq d_1\Big(\ln \frac{r_0}{|x|} +1\Big). \quad {\rm in}\ \, B_{r_2} \setminus\{0\},
$$
then for   $\epsilon>0$ 
$$u_0(x)\leq \tilde d_1 |x|^{-\epsilon} \quad {\rm in}\ \, B_{r_2} \setminus\{0\},
$$ 
which implies that
 $$
(-\Delta)^s  u_0(x)  \leq c_{24}|x|^{\theta -\epsilon p}    \quad {\rm in}\ \, B_{r_2}. \setminus\{0\}
$$
 For $\epsilon>0$  small enough,  we apply Lemma \ref{teo 2.3} with $\mu=0$ to obtain that
$$u_0(x)\leq d_1 \quad {\rm in}\ \, B_{r_2} \setminus\{0\}$$
  and we are done.

If $\tau_1<0$, applying Lemma \ref{teo 2.3} with $\mu=0$ to obtain that
$$u_0(x)\leq d_1 |x|^{\tau_1}  \quad {\rm in}\ \, B_{r_2} \setminus\{0\}$$
and $$
(-\Delta)^s  u_0(x)  \leq c_{24}|x|^{\tau_1 p+\theta}    \quad {\rm in}\ \, B_{r_2} \setminus\{0\},
$$
thus, letting $$\tau_2:=p\tau_1 +\theta +2s,$$
if $\tau_2\geq 0$,
we can get that $u_0$ is bounded at the origin  by the same proof of the case $\tau_1\geq0$ and we are done.

  If $\tau_2< 0$, iteratively, we can prove that
$$u_0(x)\leq d_1(|x|^{\tau_n}+1)\quad {\rm in}\ \, B_{r_n} \setminus\{0\},$$
where
$$\tau_n:=p\tau_{n-1}+\theta+ 2s,\quad n=1,2,\cdots .$$
 Note that
$$\tau_1-\tau_0=(p-1)\tau_0+\theta +2s>0,$$
thanks to $\tau_0>-\frac{2s}{p-1}$,
and
 \begin{eqnarray*}
\tau_n-\tau_{n-1}&=& p  (\tau_{n-1}-\tau_{n-2})\\[1.5mm]&=&p^{n-1}(\tau_1-\tau_0)
 \to+\infty\quad{\rm as}\ \, n\to+\infty,
  \end{eqnarray*}
 then there exists $n_1\in\N$ such that $\tau_{n_1}\geq 0$,
by the same argument of the case $\tau_1\geq0$, we then obtain that  $u_0$ is bounded at the origin   which contradicts  the assumption
$$\limsup_{|x|\to0^+} u_0(x) =+\infty. $$
Thus,  we conclude that
\begin{eqnarray*}
 \liminf_{|x|\to0^+}u(x) |x|^{\frac{2s+\theta}{p-1}}
  \leq  \cK_{p,\theta}   \leq \limsup_{|x|\to0^+}u(x)  |x|^{\frac{2s+\theta}{p-1}}.
\end{eqnarray*}
We complete the proof. \hfill$\Box$\medskip

Now we deal with the Sobolev's critical and supercritical case.

   \begin{proposition} \label{teo 1-sup}
Assume that  $h\in L^1_s(\R^N\setminus B_1)$,   $\theta\in(-2s,+\infty)$ and
  $$p\geq \frac{N+2s+2\theta}{N-2s}. $$

Let $u_0$ be a  positive solution of (\ref{eq 1.1}) verifying
$$\limsup_{|x|\to0^+} u_0(x) =+\infty, $$
 then
  \begin{eqnarray*}
   \liminf_{|x|\to0^+}u_0(x)  |x|  ^{\frac{2s+\theta}{p-1}}
  \leq  \cK_{p,\theta}   \leq \limsup_{|x|\to0^+}u_0(x) |x|  ^{\frac{2s+\theta}{p-1}} .
\end{eqnarray*}
 \end{proposition}
\noindent {\bf Proof. }
{\bf Part 1:}  we claim that
 \begin{equation}\label{est low}
 \liminf_{|x|\to0^+} u_0(x)|x|^{\frac{2s+\theta}{p-1}}\leq c_*^{-\frac1{p-1}} \cK_{p,\theta} .
 \end{equation}

Note that for $p\geq\frac{N+2s+2\theta}{N-2s}$,
$$-\frac{2s+\theta}{p-1}\geq \frac{2s-N}{2}.$$
Let
 $$u_p(x)= \cK_{p,\theta}  |x|^{-\frac{2s+\theta}{p-1}},$$
 where
 $$\cK_{p,\theta}=\cC_s(-\frac{2s+\theta}{p-1})^{\frac1{p-1}}.$$
 Note that
 $$(-\Delta)^s u_p= |x|^{\theta}u_p^p\quad{\rm in}\ \,\R^N\setminus\{0\},$$
which could be written as
 $$\cL_{- \cK_{p,\theta}^{p-1}}^s u_p:=(-\Delta)^s u_p- \cK_{p,\theta}^{p-1}|x|^{-2s} u_p=0\quad{\rm in}\ \, \R^N\setminus\{0\}$$
and
$$\tau_-(-\cK_{p,\theta}^{p-1})=2s-N+\frac{2s+\theta}{p-1}\quad {\rm and}\quad \tau_+(-\cK_{p,\theta}^{p-1})=-\frac{2s+\theta}{p-1}.$$

By contradiction, we assume that $u_0$ is a positive solution of (\ref{eq 1.1})  such that
$$\liminf_{|x|\to0^+}u_0(x)|x|^{\frac{2s+\theta}{p-1}}> \cK_{p,\theta}.$$
Then there exist $\epsilon_0>0$ and $r_0>0$ such that
$$|x|^\theta u_0^{p-1}(x)\geq ( \cK_{p,\theta}^{p-1}+2\epsilon_0)|x|^{-2s}\quad {\rm for\ any}\ 0<|x|<r_0.$$
Therefore, $u_0$ is a super solution of
$$\cL_{\mu_1}^s u_0\geq \epsilon_0|x|^{\theta}u_0^p \quad {\rm in}\ B_{r}\setminus\{0\},$$
 where $r\in(0,r_0)$ and
 $$\mu_1=- (\cK_{p,\theta}^{p-1}+\epsilon_0)<-  \cK_{p,\theta}^{p-1}.$$

 By \cite[Proposition 1.2]{CW},    the mapping $\mu\in(\mu_0,0) \mapsto \tau_+(\mu)$   is continuous and strictly increasing, then
  $$\tau_+( \mu_{1})<\tau_+(-  \cK_{p,\theta}^{p-1})= -\frac{2s+\theta}{p-1}\quad{\rm for}\ \ p\geq \frac{N+2s+2\theta}{N-2s}.$$
 By Lemma \ref{teo 2.2} with $\mu=\mu_1$, there exists $c>0$ such that
 $$u_0(x)\geq c_{25} |x|^{\tau_+(\mu_1)}\quad {\rm in}\ B_r\setminus\{0\}. $$

 Let $\tau_0=\tau_+(\mu_1)$.   If $\tau_0\leq -\frac{N+\theta}{p}$, where $-\frac{N+\theta}{p}> \frac{2s-N}2$ for  $p\geq \frac{N+2s+2\theta}{N-2s}$.  Then $|x|^{\theta} u_0^p\not\in L^1(B_1)$, a contradiction arises from Theorem \ref{theorem-C} part $(iii)$
 with $\mu=0$. Then we are done.

If $\tau_0\in\big(-\frac{N+\theta}{p}, -\frac{2s+\theta}{p-1} \big)$,  then
 $$
(-\Delta)^s  u_0(x)  \geq    d_0^p|x|^{p\tau_0+\theta}=  d_0^p |x|^{\tau_1-2s}\quad {\rm in}\ \, B_{r_0} \setminus\{0\},
$$
where
$$\tau_1:=p\tau_0+\theta+2s.$$

If $\tau_1 \leq -\frac{N+\theta}{p}$, we are done.

If $\tau_1\in\big(-\frac{N+\theta}{p},-\frac{2s+\theta}{p-1} \big)$,  by Lemma \ref{teo 2.2} with $\mu=0$, we have that
$$u_0(x)\leq d_1|x|^{\tau_1}\quad {\rm in}\ \, B_{r_0} \setminus\{0\}.$$
Iteratively, we recall that
$$\tau_j:=p\tau_{j-1}+\theta+2s,\quad j=1,2,\cdots .$$

If $\tau_{j} \leq -\frac{N+\theta}{p}$ we are done, otherwise
it follows from Lemma \ref{teo 2.2}  that
$$u_0(x)\geq d_{j+1} |x|^{\tau_{j+1}},$$
where
$$\tau_{j+1}=p\tau_j+\theta+2s<\tau_j.$$

Thanks to  $\tau_1-\tau_0 =(p-1)\tau_0+\theta+2s<0$, we have that
$$
\tau_j-\tau_{j-1} = p(\tau_{j-1}-\tau_{j-2})=p^{j-1} (\tau_1-\tau_0)\to-\infty\ \, {\rm as}\ j\to+\infty,
$$
 then there exists $j_0\in\N$ such that
$$\tau_{j_0}>-\frac{N+\theta}{p}\quad {\rm and}\quad \tau_{j_0+1}\leq  -\frac{N+\theta}{p}.$$
For this,  we have  that
\begin{eqnarray*}
(-\Delta)^s  u_0(x) \geq      d_{j_0}^p  |x|^{\tau_{j_0+1}-2s}  \quad{\rm in}\ \ B_{r_0}\setminus\{0\},
 \end{eqnarray*}
 then $u_0(x)\geq d_{j_0+1} |x|^{\tau_{j_0+1}},$
 and  a contradiction arises. \smallskip\smallskip

 {\bf Part 2: } we claim
$$ \limsup_{|x|\to0^+} u_0(x)|x|^{\frac{2s+\theta}{p-1}}\geq \cK_{p,\theta}.$$

By contradiction, we assume that $u_0$ is a positive super solution of (\ref{eq 1.1}) in $\Omega\setminus\{0\}$
such that
$$\limsup_{|x|\to+\infty}u_0(x)|x|^{\frac{2s +\theta}{p-1}}< \cK_{p,\theta}.$$
Then there exist $\epsilon_1\in(0,\frac12 \cK_{p,\theta}^{p-1})$ and $r>0$ such that
$$|x|^\theta u_0^{p-1}(x)\leq \big( \cK_{p,\theta}^{p-1}-\epsilon_1\big)|x|^{-2s}\quad {\rm for\ any}\ 0<|x|<r.$$
Therefore, $u_0$ is a super solution of
$$\cL_{\mu_{2}}^s u_0\leq 0 \quad {\rm in}\ B_{r_0}\setminus\{0\},$$
 where
 $$\mu_{2}:=- \cK_{p,\theta}^{p-1}+ \epsilon_1>-   \cK_{p,\theta}^{p-1}. $$

  Notice that for $p\geq \frac{N+2\theta+2s}{N-2s}$, we have that
 $$\tau_+(-\cK_{p,\theta}^{p-1})=-\frac{2s+\theta}{p-1}\geq \frac{2s-N}{2}>\tau_-\big(- \cK_{p,\theta}^{p-1}\big).$$
 By the  strictly increasing monotonicity, we have that
 $\tau_+(\mu_{2})>-\frac{2s+\theta}{p-1},$
 thank to  $\mu_{2}>-\cK_{p,\theta}^{p-1}$.  By Lemma \ref{teo 2.3} with $\mu=\mu_2$, we see that
 $$u_0(x)\leq c_{26}|x|^{\tau_+(\mu_{2})}\quad {\rm in}\ B_{r_0}\setminus\{0\}. $$

Let $\tau_0=\tau_+(\mu_{2})>-\frac{2s+\theta}{p-1}$.   Since $\tau_0\in\big(-\frac{2s+\theta}{p-1}, 0\big)$,  then
 $$
(-\Delta)^s  u_0(x)  \leq    c_{27} |x|^{p\tau_0+\theta}= c_{27}  |x|^{\tau_1-2s}\quad {\rm in}\ \, B_{r_0} \setminus\{0\},
$$
where
$$\tau_1:=p\tau_0+\theta+2s.$$
if $\tau_1> 0$,
 $$u_0(x)\leq d_1 \quad {\rm in}\ \, B_{r_2} \setminus\{0\},$$
 which contradicts  that $u_0$ is non-removable at the origin. Then we are done.

If $\tau_1=0$, applying Lemma \ref{teo 2.3} with $\mu=0$ to obtain that
$$u_0(x)\leq d_1\Big(\ln \frac{r_0}{|x|} +1\Big)\quad {\rm in}\ \, B_{r_2} \setminus\{0\},
$$
then for $\epsilon>0$ small,  
$$u_0(x)\leq \tilde d_1 |x|^{-\epsilon p} \quad {\rm in}\ \, B_{r_2} \setminus\{0\},
$$
which implies that
 $$
(-\Delta)^s  u_0(x)  \leq c_{24}|x|^{\theta -\epsilon p}    \quad {\rm in}\ \, B_{r_2} \setminus\{0\}.
$$
 Then we apply Lemma \ref{teo 2.3} with $\mu=0$ to obtain that
$$u_0(x)\leq d_1 \quad {\rm in}\ \, B_{r_2} \setminus\{0\}$$
for $\epsilon>0$ small enough and we are done.

If $\tau_1< 0$,  by Lemma \ref{teo 2.3}, we have that
$$u_0(x)\leq d_1|x|^{\tau_1}\quad {\rm in}\ \, B_{r_0} \setminus\{0\}.$$
Iteratively, we reset that
$$\tau_j:=p\tau_{j-1}+\theta+2s,\quad j=1,2,\cdots .$$
Note that
$$\tau_1-\tau_0=(p-1)\tau_0+\theta+2s >0$$
thanks to $\tau_0>-\frac{2s+\theta}{p-1}$.\smallskip

If $\tau_{j}p   +\theta +2s\geq 0$ we are done by the same argument as $\tau_1\geq0$.
If $\tau_{j}p   +\theta +2s< 0$
it follows by Lemma \ref{teo 2.2}  that
$$u_0(x)\geq d_{j+1} |x|^{\tau_{j+1}},$$
where
$$\tau_{j+1}=p\tau_j+2s+\theta<\tau_j.$$

Note that
$$
\tau_j-\tau_{j-1} = p(\tau_{j-1}-\tau_{j-2})=p^{j-1} (\tau_1-\tau_0)\to+\infty\ \, {\rm as}\ j\to+\infty,
$$
 thus, there exists $j_0\in\N$ such that
$$\tau_{j_0}< 0\quad {\rm and}\quad \tau_{j_0+1}\geq 0,$$
and  $u_0$ is bounded at the origin.
Thus  a contradiction comes from the assumption that
 $$\limsup_{|x|\to0^+} u_0(x) =+\infty. $$
 We complete the proof.   \hfill$\Box$\medskip

 \subsection{The isolated singularities}

   \noindent{\bf Proof of Theorem \ref{teo 0}. }  Let  $u_0$ be a nonnegative solution of (\ref{eq 1.1}).  For
  $$\theta=0 \quad \&\quad  p\in \Big(\frac{N }{N-2s}, \frac{N+2s }{N-2s}\Big),$$
Proposition \ref{pr 2.1} gives an upper bound
$$
 u_0(x)\leq c_{16} |x|^{-\frac{2s}{p-1}} ,\quad\forall\, x\in B_1\setminus\{0\}.
$$
Then $u_0$ verifies Harnack inequality (\ref{harn in}) by Proposition \ref{pr Harnack}.
We now apply Proposition \ref{teo 1-sub} with
$p\in \big(\frac{N}{N-2s}, \frac{N+2s}{N-2s}\big)$,
  we obtain that
 \begin{eqnarray*}
\frac{\cK_{p}}{C_0} \leq \liminf_{|x|\to0^+}u(x) |x|^{\frac{2s}{p-1}}
  \leq  \cK_{p}   \leq \limsup_{|x|\to0^+}u(x) |x|^{\frac{2s}{p-1}}\leq \cK_{p} C_0.
\end{eqnarray*}
 We complete the proof.    \hfill$\Box$ \medskip

\noindent{\bf Proof of  Theorem \ref{teo 1}. } When $$\theta=0 \quad\&\quad  p\geq \frac{N+2s}{N-2s}$$
 or
 $$\theta\in(-2s,0) \quad\&\quad p>\frac{N+\theta}{N-2s},$$
 and  $h$  is radially symmetric, decreasing with respect to $|x|$,  then the positive solution $u_0$ is radially symmetric and decreasing with respect to $r=|x|$, provided that $u_0(x)\geq h(1)$,
then Proposition \ref{pr 2.1-sup} provides
   an upper bound
$$
 u_0(x)\leq c_{22} |x|^{-\frac{2s+\theta}{p-1}} ,\quad\forall\, x\in B_1\setminus\{0\}.
$$
Then our conclusion of $u_0$ follows by Proposition \ref{pr Harnack},  Proposition \ref{teo 1-sub}  and Proposition \ref{teo 1-sup}.
      \hfill$\Box$ \medskip

\noindent{\bf Proof of  Theorem \ref{pr 1}. }  For $\theta\in(-2s,+\infty)$, Corollary \ref{cr 3.1-1}
gives an upper bound
$$
 u_0(x)\leq c_{22} |x|^{-\frac{2s+\theta}{p-1}} ,\quad\forall\, x\in B_1\setminus\{0\}.
$$
Then our conclusion of $u_0$ follows by   Proposition \ref{teo 1-sub}  and Proposition \ref{teo 1-sup}.
      \hfill$\Box$ \medskip

\noindent{\bf Proof of  Theorem \ref{teo 2}. } Let $v_0$ be a positive solution of (\ref{eq 1.1-lap}). We claim that
\begin{equation}\label{ess-pp}
v_0(x)\leq c_{28} |x|^{-\frac{2+\theta}{p-1}}\quad {\rm for}\ \, 0<|x|< 1/2
\end{equation}
for some $c_{28}>0$.

We first show that for $c_{29}>0$
\begin{equation}\label{ess-p}
\int_{B_r} |x|^{\theta} v_0^p  dx\leq c_{29} r^{N-\frac{\theta +2p}{p-1}}, \quad \forall \ r\in(0,\frac{1}{2}).
\end{equation}

In fact,  let $r\in(0,1/2)$ and
 $$v=v_0-\tilde v_0,$$
 where
 $\tilde v_0$ is the  harmonic extension of $v_0$ in $B_r$, i.e. the solution of
$$
 \left\{\arraycolsep=1pt
\begin{array}{lll}
 -\Delta   u= 0\quad\ \
   &{\rm in}\ \,  B_r,\\[2mm]
\quad \ \ u= v_0\quad     &{\rm   on}\ \,  \partial B_r.
  \end{array}
 \right.
$$
 Then we have that
$$
 \left\{\arraycolsep=1pt
\begin{array}{lll}
 -\Delta   v=|x|^{\theta} (v+\tilde v_0)^p\quad\ \
   &{\rm in}\ \,  B_r\setminus\{0\},\\[2mm]
\quad \ \ v= 0\quad \ \  &{\rm   on}\ \,  \partial B_r
  \end{array}
 \right.
$$
and \cite{L} shows that in the Serrin's supercritical  since $p>\frac{N+\theta}{N-2}$,
$$
  \int_{B_r}v(-\Delta)  \eta dx=\int_{B_r}|x|^{\theta} (v+\tilde  v_0)^p  \eta dx,\quad \eta\in C^{1.1}_0(B_r),
$$
  which implies that for $\eta\in C^{1.1}_0(B_r)$
  \begin{eqnarray*}
  \int_{B_r}v_0(-\Delta)  \eta dx&=&\int_{B_r} |x|^{\theta} v_0^p\eta dx+\int_{B_r} \tilde v_0 (-\Delta)  \eta dx
  \\&=&\int_{B_r} |x|^{\theta} v_0^p\eta dx-\int_{\partial B_r}  v_0 \nabla\eta\cdot x d\sigma(x),
\end{eqnarray*}
since  the normal vector pointing outside of $B_1$ is $x\in \partial B_1$.

Now we take $\eta=\xi_r:=\xi_{1}(r^{-1}\cdot)$,
where $\xi_1$ is the first eigenvalue of $-\Delta$ in $B_1$ with the zero Dirichlet boundary condition. Note that  $\xi_1$ is radially symmetric and decreasing with respect to $|x|$
and then $\xi_1'(1)<0$.

Then we have that
 $$
  \lambda_1 r^{-2}\int_{B_r}v_0 \xi_r dx=\int_{B_r} |x|^{\theta} v_0^p \xi_r dx-r^{-1}\xi_1'(1) \int_{\partial B_r}  v_0 d\sigma(x)>\int_{B_r} |x|^{\theta} u_0^p \xi_r dx.
 $$

  We observe that $-\frac{\theta}{p-1}>-N$ by our assumption that $p>\frac{N+\theta}{N-2}$ and
  \begin{eqnarray*}
\int_{B_r} |x|^{\theta} v_0^p \xi_r dx &<&  \lambda_1 r^{-2}\int_{B_r}v_0\xi_r dx
\\[1.5mm]&\leq &\lambda_1 r^{-2}\Big(\int_{B_r}  |x|^{\theta} v_0^p  \xi_r dx\Big)^{\frac1p}\Big(\int_{B_r}\xi_r|x|^{-\frac{\theta}{p-1}}  dx\Big)^{1-\frac1p}
\\[1.5mm]&\leq & \|\xi_1\|_{L^\infty(B_1)}^{1-\frac1p}\lambda_1 r^{(N-\frac{\theta}{p-1})(1-\frac1p)-2}\Big(\int_{B_r}|x|^{\theta}v_0^p \xi_r dx\Big)^{\frac1p},
\end{eqnarray*}
which implies that
$$
\int_{B_r} |x|^{\theta} v_0^p\xi_r dx < \|\xi_1\|_{L^\infty(B_1)} \lambda_1^{\frac{p}{p-1}} r^{N -\frac{\theta}{p-1}-\frac{ 2  p}{p-1}}.
$$
Notice   that
$$\xi_r(x)=\xi_1(r^{-1}x)\geq \min_{z\in\bar B_{\frac12}}\xi_1(z)\quad {\rm in }\ B_{\frac r2},$$
then
 $$
\int_{B_\frac r2}  |x|^{\theta} v_0^p  dx< c_{30}\int_{B_r}  |x|^{\theta} v_0^p\xi_r dx < c_{30}\|\xi_1\|_{L^\infty(B_1)} \lambda_1^{\frac{p}{p-1}} r^{N -\frac{\theta}{p-1}-\frac{ 2  p}{p-1}}.
 $$
 Replacing $r$ by $\frac r2$, we obtain (\ref{ess-p}). \smallskip

{\it Part $(i)$: } Since $v_0$ verifies (\ref{lp-2--}), thus, taking $|x|=r$,  for any $y\in B_r\setminus B_{\frac{r}2}$, there holds that
 $$u_0(y)\geq \frac1{C_1}u_0(x). $$
  Now we deduce from  (\ref{ess-p})  that for any given $|x|=r$,
   \begin{eqnarray*}
C_0^{-p}\Big(1-(\frac12)^N\Big) |\partial B_1| \, v_0^p(x) r^{N-\theta}\leq  \int_{B_r\setminus B_\frac r2} |y|^\theta v_0^p(y)  dy<c_{14} r^{N-\frac{\theta}{p-1}- \frac{ 2s p}{p-1}},
 \end{eqnarray*}
 which implies that for some $c_{31}\geq 1$
 $$\frac1{c_{31}} |x|^{ -\frac{2+\theta}{p-1}}\leq v_0 (x)\leq c_{13} |x|^{ -\frac{2+\theta}{p-1}}\quad {\rm for}\ \,  0<|x|<\frac12. $$
By similar argument in Proposition \ref{teo 1-sup} (some basic estimates could see \cite{CQZ}), we obtain (\ref{ss-s}).
\smallskip

{\it Part $(ii)$: } Let $u=v_0-h\geq0$, under our assumptions that $h$ is a nonnegative constant and
$v_0\geq h$ in $B_1$,
then $u$ is a classical solution of
\begin{equation}\label{eq 1.1-0lap}
\left\{
\begin{array}{lll}
 -\Delta  u=  |x|^{\theta} \big(u+ h\big)^{p}  \quad &{\rm in}\ \, B_1 \setminus\{0\},
\\[2mm]
\quad \ \ \, u=0&{\rm on}\ \,  \partial B_1.
\end{array}\right.
 \end{equation}
 For $\theta<0$, the moving plane  works for the problem (\ref{eq 1.1-0lap}), see \cite{BN,GNN} and we get that
 $u$ is radially symmetric and decreasing with respect to $|x|$, so is $v_0$. Then   (\ref{ess-p}) implies (\ref{ess-pp}). \smallskip

Under the upper bound (\ref{ess-pp}), we have  the Harnack inequalities,  for all $r\in(0,\frac12)$
$$
 \sup_{x\in B_{2r}\setminus B_r} v_0(x)\leq C_1 \inf_{x\in B_{2r}\setminus B_r} v_0(x)
$$
 for some $C_1>0$. From the Harnack inequality  in \cite{S}, we have that for some constant $c_{32}\geq1$
  $$\frac1{c_{32}}|x|^{-\frac{2+\theta}{p-1}}\leq v_0(x)\leq c_{32}|x|^{-\frac{2+\theta}{p-1}}\quad {\rm for}\ \,  0<|x|<\frac12.$$
Therefore, we obtain (\ref{ss-s}).   \hfill$\Box$

      \setcounter{equation}{0}
    \section{Application to the decay at infinity}
    In this section,  we apply the classification of isolated singularity to obtain the decay at infinity for the
    positive solution of
\begin{equation}\label{eq 1.1-ex}
\left\{
\begin{array}{lll}
(-\Delta)^s u=  |x|^{\tilde \theta} u^{p}  \quad &{\rm in}\ \,  \R^N\setminus \bar B_1,
\\[2mm]
\qquad \ \ u=h&{\rm in}\ \,   \bar B_1.
\end{array}\right.
 \end{equation}

 \begin{theorem} \label{teo 1-ex}
  Assume that     
  \begin{equation}\label{con ss-0-e}
  p\in \Big(\frac{N+\tilde \theta}{N-2s},\, \frac{N+2s+\tilde \theta}{N-2s}\Big)\qquad   {\rm for }  \quad \tilde \theta>-2s
\end{equation}
  or
  \begin{equation}\label{con ss+0-e}
p=\frac{N+2s+\tilde \theta}{N-2s}  \qquad {\rm for }  \quad \tilde \theta\geq 0, 
\end{equation}
$h\in L^1( \bar B_1)$ is a nonnegative and radially symmetric function such that $h(1)\geq0$,    $|x|^{N-2s}h(x)$ is increasing with respect to $|x|$ in $B_1\setminus\{0\}$.

Let
$u_0$ be a  positive solution of (\ref{eq 1.1-ex})   such that
 \begin{equation}\label{con h-c}
  u_0(x)\geq  h(1)|x|^{2s-N}\quad{\rm for}\ \, \forall\, x\in \R^N\setminus B_1,
  \end{equation}
then  $u_0$ is radially symmetric,
   either  for some $k\geq h(1)$
  \begin{eqnarray*}
 \lim_{|x|\to+\infty}u_0(x)  |x|^{N-2s}=k
\end{eqnarray*}
or
  \begin{eqnarray*}
\frac{\cK_{p,\theta^*}}{C_0} \leq  \liminf_{|x|\to+\infty}u_0(x)  |x|^{ \frac{2s+\tilde\theta }{p-1}}
  \leq  \cK_{p,\theta^*}   \leq \limsup_{|x|\to+\infty}u_0(x)  |x|^{\frac{2s+\tilde\theta }{p-1}}  \leq  C_0\cK_{p,\theta^*},
\end{eqnarray*}
where
$$\theta^*=p(N-2s)-N-2s-\tilde \theta.$$

 \end{theorem}
 \noindent{\bf Proof. } We use the  {\it Kelvin transformation: }
     let $u_f$ be a  solution of
     \begin{equation}\label{eq 2.1-ex}
\left\{\arraycolsep=1pt
\begin{array}{lll}
(-\Delta)^s u=  f  \quad\       &{\rm in}\ \,   \cO , \\[2mm]
\qquad \ \  u=h\quad  \   & {\rm in}\ \    \R^N\setminus \cO,
  \end{array}
 \right.
\end{equation}
  where $\cO$ is an open set of $\R^N$ and $f$  is locally  H\"older continuous  in $\cO$.
   Denote
$$
   u^\sharp(x)=|x|^{2s-N}u_f(\frac{x}{|x|^2})\quad{\rm for}\  x\in \R^N\setminus \{0\}
$$
and
$$\cO^\sharp=\Big\{x\in \R^N: \frac{x}{|x|^2}\in \cO\Big\},$$
then a  direct computation shows that
$$
 (-\Delta)^s u^\sharp(x)    =  |x|^{-2s-N} \big((-\Delta)^s u_f\big)\left(\frac{x}{|x|^2}\right)\quad{\rm for}\  x\in \cO^\sharp.
$$
As a conclusion,  $u^\sharp$ is a solution of
 $$
\left\{\arraycolsep=1pt
\begin{array}{lll}
(-\Delta)^s  u^\sharp = |x|^{-2s-N}  f^* \quad\       &{\rm in}\ \,    \cO^\sharp, \\[2mm]
\qquad \ \  u^\sharp = |x|^{2s-N}  h^* \quad  \   & {\rm in}\ \    \R^N \setminus  \cO,
  \end{array}
 \right.
$$
where $  f^*(x)=f(\frac{x}{|x|^2})$ and $h^*(x)=h(\frac{x}{|x|^2})$.

Let $u_0$ be a positive solution of (\ref{eq 1.1-ex}).   Taking $f(x)=|x|^{\tilde \theta}u_0^{p}(x)$ and $\cO=\R^N\setminus B_1$
in (\ref{eq 2.1-ex}),    the function
$$u^\sharp(x):=|x|^{2s-N}u_0(\frac{x}{|x|^2})$$
 is a positive solution of
 $$
\left\{\arraycolsep=1pt
\begin{array}{lll}
(-\Delta)^s  u^\sharp =|x|^{\theta^*}  (u^\sharp)^p \quad\       &{\rm in}\ \,   B_1\setminus \{0\} ,
\\[2mm]
\qquad \ \  u^\sharp=|x|^{2s-N} h^*\quad  \   & {\rm in}\ \     \R^N\setminus B_1,
  \end{array}
 \right.
$$
where
$$
 |x|^{2s-N} h^*\in L^1(\R^N\setminus B_1, |x|^{-N-2s}dx), $$
 which is equivalent to
 $h\in  L^1(B_1)$.

Since
$$
 |x|^{2s-N} h^*(x)=|y|^{N-2s}h(y)\quad{with}\ \, x=\frac{y}{|y|^2}, $$
 then $x\in  B_1\setminus\{0\} \mapsto |x|^{2s-N} h^*(x)$ is radially symmetric and increasing
 with respect to $|x|$.
 Note that
 $$   u_0(x)\geq h(1)|x|^{2s-N} \quad {\rm for}\ \,  \forall x\in \R^N\setminus B_1\quad\ \Longleftrightarrow \quad\    u^\sharp(y)\geq h(1)\quad {\rm for}\ \,  y\in B_1\setminus\{0\}.$$
To obtain the asymptotic behaviors,  we would apply Theorem \ref{teo 1}
 with
  $\theta^*=0$ and $p\geq \frac{N+2s}{N-2s}$  by (\ref{con ss+0-e});
  or with
 $\theta^*\in(-2s,0)$ and $p\geq \frac{N+\theta^*}{N-2s}$  by  (\ref{con ss-0-e}).
 As a consequence, we have that $u^\sharp$ is radially symmetric, decreasing with respect to $|y|$, and is either removable at the origin or
verifies
 \begin{eqnarray*}
\frac{\cK_{p,\theta^*}}{C_0} \leq  \liminf_{|x|\to0^+}u^\sharp(x)  |x|  ^{\frac{2s+\theta^*}{p-1}}
  \leq  \cK_{p,\theta^*}   \leq \limsup_{|x|\to0^+}u^\sharp(x) |x|  ^{\frac{2s+\theta^*}{p-1}}  \leq  C_0\cK_{p,\theta^*},
\end{eqnarray*}
 which  imply that $u_0$ is radially symmetric, 
 either
 \begin{eqnarray*}
 \lim_{|x|\to+\infty}u_0(x)  |x|^{N-2s}=u^\sharp(0)
\end{eqnarray*}
or
  \begin{eqnarray*}
\frac{\cK_{p,\theta^*}}{C_0}  \leq   \liminf_{|x|\to+\infty}u_0(x)  |x|^{N- \frac{2sp+\theta^*}{p-1}}
 &\leq&  \cK_{p,\theta^*}
   \\[2mm] & \leq& \limsup_{|x|\to+\infty}u_0(x)  |x|^{N-\frac{2sp+\theta^*}{p-1}}  \leq  C_0\cK_{p,\theta^*}.
\end{eqnarray*}
We complete the proof. \hfill$\Box$\medskip

\begin{theorem} \label{teo 0-ex}
Assume that  $h\in   L^1(B_1)$
  and
   \begin{equation}\label{con ss+0-e1}
p=\frac{N+2s+\tilde \theta}{N-2s}  \qquad {\rm for }  \quad \tilde \theta\in(-2s,\, 0).
\end{equation}

  Let  $u_0$ be a  positive solution of (\ref{eq 1.1-ex}) and $C_0\geq 1$ be the constant from the Harnack inequality in Theorem \ref{teo 0} with $\theta=0$. 

  Then
   either  for some $k>0$
  \begin{eqnarray*}
 \lim_{|x|\to+\infty}u_0(x)  |x|^{N-2s}=k
\end{eqnarray*}
or
  \begin{eqnarray*}
\frac{\cK_{p}}{C_0} \leq  \liminf_{|x|\to+\infty}u_0(x)  |x|^{ \frac{2s+\tilde\theta }{p-1}}
  \leq  \cK_{p}   \leq \limsup_{|x|\to+\infty}u_0(x)  |x|^{\frac{2s+\tilde\theta }{p-1}}  \leq  C_0\cK_{p}.
\end{eqnarray*}
 \end{theorem}
 \noindent{\bf Proof. }
  From the Kelvin transformation,
$$u^\sharp(x):=|x|^{2s-N}u_0(\frac{x}{|x|^2})$$
 is a positive solution of
 $$
\left\{\arraycolsep=1pt
\begin{array}{lll}
(-\Delta)^s  u^\sharp =   (u^\sharp)^p \quad\       &{\rm in}\ \,   B_1\setminus \{0\} ,
\\[2mm]
\qquad \ \  u^\sharp=|x|^{2s-N} h^*\quad  \   & {\rm in}\ \     \R^N\setminus B_1,
  \end{array}
 \right.
$$
by our assumption (\ref{con ss+0-e1}) that  $\tilde \theta=p(N-2s)-N-2s$.

From Theorem \ref{teo 0}, $u^\sharp$ is either removable at the origin or it
verifies
 \begin{eqnarray*}
\frac{\cK_{p}}{C_0} \leq  \liminf_{|x|\to0^+}u^\sharp(x)  |x|  ^{\frac{2s }{p-1}}
  \leq  \cK_{p}   \leq \limsup_{|x|\to0^+}u^\sharp(x) |x|  ^{\frac{2s}{p-1}}  \leq  C_0\cK_{p},
\end{eqnarray*}
then we have that
 either
 \begin{eqnarray*}
 \lim_{|x|\to+\infty}u_0(x)  |x|^{N-2s}=u^\sharp(0)
\end{eqnarray*}
or
  \begin{eqnarray*}
\frac{\cK_{p}}{C_0}  \leq   \liminf_{|x|\to+\infty}u_0(x)  |x|^{N- \frac{2sp }{p-1}}
  &\leq&  \cK_{p}
   \\[2mm] & \leq& \limsup_{|x|\to+\infty}u_0(x)  |x|^{N-\frac{2sp }{p-1}}  \leq  C_0\cK_{p},
\end{eqnarray*}
where we use $\theta^*=p(N-2s)-N-2s-\tilde \theta=0.$
We complete the proof. \hfill$\Box$\medskip

\bigskip
 \bigskip

   \noindent{\bf \small Acknowledgements:} {\footnotesize Foundation of China, No. 12071189, by Jiangxi Province Science Fund for Distinguished Young Scholars, No. 20212ACB211005, and by the Jiangxi Provincial Natural Science Foundation, No. 20202ACBL201001, by the Science and Technology Research Project of Jiangxi Provincial Department of Education, No. GJJ200307, GJJ200325.  F. Zhou is supported
by Science and Technology Commission of Shanghai Municipality (STCSM), Grant No. 18dz2271000
and also supported by NSFC (No. 11431005). }

  \end{document}